\documentclass[11pt,reqno]{amsart}

\usepackage{euscript}
\usepackage{pb-diagram}
\usepackage{lamsarrow}
\usepackage{pb-lams}
\usepackage{epic}
\usepackage{epsfig}
\usepackage{amssymb}

\theoremstyle{plain}

\theoremstyle{definition}

\newcommand{\Gm}{\G^{min}(\xfn,z)}
\newcommand{\Gc}{\G^{can}(\xfn,z)}
\newcommand{\coker}{\mbox{coker}\,}
\newcommand{\bfc}{{\mathbb C}}
\newcommand{\bfz}{{\mathbb Z}}
\newcommand{\bfq}{{\mathbb Q}}

\newcommand{\calst}{{\mathcal  St}}
\newcommand{\calg}{{\mathcal  G}}
\newcommand{\cald}{{\mathcal  D}}
\newcommand{\call}{{\mathcal  L}}
\newcommand{\calw}{{\mathcal  W}}
\newcommand{\calv}{{\mathcal  V}}
\newcommand{\cala}{{\mathcal  A}}
\newcommand{\calr}{{\mathcal  R}}
\newcommand{\G}{\Gamma}
\newcommand{\xfn}{X_{f,n}\,}
\newcommand{\bc}{{\mathbb C}}
\newcommand{\bz}{{\mathbb Z}}

\newcommand{\q}{{\mathbb Q}}

\newcommand{\finv}{f^{-1}(0)} 
\newcommand{\m}{\mathtt{m}}   
\newcommand{\g}{\mathtt{g}}

\begin{document}      

\title{The link of $\{f(x,y)+z^n=0\}$ and Zariski's conjecture}

\author{Robert Mendris and Andr\'as N\'emethi}
\address{Department of Mathematics\\Ohio State University\\Columbus, OH 43210}
\email{mendris@math.ohio-state.edu}
\email{nemethi@math.ohio-state.edu}
\thanks{The second author is partially supported by NSF grant
DMS-0071820.}
\begin{abstract} We consider suspension hypersurface singularities of 
type $g=f(x,y)+z^n$, where $f$ is an irreducible plane curve singularity.
For such germs, we prove that the link of $g$ determines  completely the 
Newton pairs of $f$ and the integer $n$ except for  two pathological cases,
which can be completely described. Even in the pathological  cases, the link 
and the Milnor number of $g$  determine uniquely the Newton pairs of $f$
and $n$. In particular, for such $g$,
we verify Zariski's conjecture about the multiplicity. 
The result also  supports the following conjecture formulated in the 
paper. 
If the link of an isolated hypersurface singularity is a rational homology 
3-sphere then it determines the embedded topological type, the equivariant
Hodge numbers and the multiplicity of the singularity. The conjecture is 
verified for weighted homogeneous singularities too. 

\end{abstract}

\maketitle
\pagestyle{myheadings}
\markboth{{\normalsize Robert Mendris and Andr\'as N\'emethi}}{
{\normalsize On the link of $\{f(x,y)+z^n=0\}$}}

\section{Introduction.}

In the last decades an intense research effort has been
 concentrated on the following
problem: what kind of analytic invariants  or smoothing invariants (if they 
exist) can be determined 
from the topology of a normal surface singularity $(X,x)$.

Some of the results became already classical. E.g., Mumford's result which
states that $(X,x)$ is smooth if and only if the fundamental group of the 
link $L_X$ is trivial \cite{Mumford}. Or its generalization by Neumann
\cite{NePl} which claims that the oriented homeomorphism type of the link
contains the same information as the resolution graph of $(X,x)$. 
Or Artin's computations of the multiplicity and the embedded dimension
of rational singularities \cite{Artin,Artin2}; and their  generalizations 
by Laufer \cite{Lauferme} for minimally elliptic singularities, and by
S. S.-T. Yau \cite{Yau} for some elliptic  singularities. 

In general, these questions are very difficult, even if we restrict ourselves 
to some special families, e.g. to complete intersections or
hypersurface singularities;  and even if we permit ourselves to use, instead of
the topology of $(X,x)$,  a  richer topological information, e.g.
in the case of hypersurfaces the  {\em embedded topological type}. 
For example, Zariski conjectured three decades ago that 
the embedded topological type of an isolated hypersurface singularity
  determines its multiplicity \cite{ZarConj}. 
This was verified till now only for quasi-homogeneous singularities 
\cite{Greuel,Osh,XY,YauConj} (and some other sporadic cases).

In a different direction, the conjecture of Neumann and Wahl \cite{NW}
about a possible connection between the Casson invariant of the link
(provided that it is an integral homology sphere) and the signature 
of the Milnor fiber opened new windows for the theory. 

Recently,  the subject revives even with a larger intensity, see e.g. 
\cite{NInv,NW2,NW3,NN,NNII,NNIII}, their introductions and listed  
references. Basically, these articles claim that
if the link of a $\q$-Gorenstein singularity  is a rational homology
sphere, then it codifies an extremely rich analytic information about the 
singularity. 

The present article is in the spirit of the above efforts. We will consider 
the family of suspension hypersurface  singularities of type $f(x,y)+z^n$, 
where $f$ is an irreducible plane curve singularity. 

For this family, we not only answer positively both main conjectures, (namely 
the Zariski's conjecture, and also the possibility to recover the main
analytic and smoothing invariants from the link), but we succeed to obtain 
much sharper statements. 

The main result of the paper is the following (cf. \ref{main}):

\vspace{2mm}

\noindent {\bf Theorem 1.} {\em Let $f:(\bc^2,0)\to(\bc,0)$ 
be an irreducible plane curve 
singularity with Newton pairs $\{(p_i,q_i)\}_{i=1}^s$ and let $n\geq 2$ be an
integer.  Let $L_X$ be the link of 
the hypersurface singularity $(X,0)=(\{f(x,y)+z^n=0\},0)$.
Then, except for two pathological cases S1  and S2 (which are described
completely in \ref{S1} and 
\ref{S2},  and can be characterized perfectly in terms of $L_X$), 
from the link $L_X$ one can recover completely the Newton pairs
of $f$ and the integer $n$ (provided that  we disregard the 
``$z$-axis ambiguity'', cf. \ref{zamb}.) 
In both exceptional cases the 
links have non-trivial first Betti numbers. In particular, 
the above statement holds  without any 
exception provided that the link is a rational homology sphere.

On the other hand, in the cases S1 and S2, the link together with 
the Milnor number of the hypersurface singularity $f+z^n$ determine
completely the Newton pairs of $f$ and the integer $n$ (cf.
the two addendums in \ref{S1} and \ref{S2}). }

\vspace{2mm}

Here some remarks are in order.

\vspace{2mm}

(1) $L_X$ determines the  number $s$ of
Newton pairs of $f$ in all the cases. 

The exceptional case S1 appears when $s=1$, and the corresponding 
singularities have the equisingular type of some 
special Brieskorn singularities. This case can be easily classified.

The exceptional case S2 appears for $s=2$ with some other strong additional 
restrictions. In this case any link $L_X$ can be realized by at most 
two possible pairs $(f,n)$. This case again is completely  clarified.

In all other cases, e.g. when $s\geq 3$, the theorem assures uniqueness.
This is slightly surprising. At the beginning of our study, we expected here
more and more complicated special families providing interesting coincidences
for their links. But, it turns out that this is not the case: if the plumbing
graph of the link (or equivalently, the resolution graph of $(X,x)$)
 has more and more complicated structure, then it becomes more 
``over-determined'', and it  leaves no room for any ambiguity
for $f$ and $n$. 

In fact, in order to reach our goal, it was sufficient to consider a rather
limited 
information about this graph: the determinants of its  maximal strings
and the determinants of some subgraphs with only one rupture vertex. 

Except for the two pathological cases (when the graphs are really very simple),
in all other cases already 
these determinants determine all the Newton pairs and $n$.

(2) In general, it is very difficult to characterize those 
resolution graphs (or links) which can be realized by, say, hypersurface 
singularities, or complete intersections, or by any family of germs 
defined by some analytic property. 

Our proof gives a complete characterization of those graphs
which can be realized as resolution graphs of some  $\{f+z^n=0\}$ for 
some irreducible $f$. Indeed, the proof is a precise recipe how one can 
recover the Newton pairs  of $f$ and $n$. If one runs this algorithm (the steps
of the proof of \ref{main})  for
an arbitrary minimal resolution graph, and at some point it fails, then the
graph definitely is not of this type. If the algorithm goes through and
provides some candidates for the 
Newton pairs of $f$  and for $n$, then one has to
compute the minimal resolution graph of $f+z^n=0$ (using e.g.
\ref{6.3}) and compare it  with the 
initial graph. If they are the same, then the answer is yes. But it can happen
that these two graphs  are not the same 
(since our algorithm is based on a very 
limited number of determinants: these determinants  can be the same even if 
the  graphs are not). 

(3) For the case when $f$ is arbitrary (i.e. reduced),  but $n=2$, 
 and the link is rational homology sphere, Laufer established the uniqueness 
in \cite{Lauferntwo}. 

In \cite{Pi1,Pi2},  A. Pichon proved by a different method
that any link can be realized in a {\em finitely many} ways as the link of 
$f+z^n=0$ (in her case $f$ is reduced too). We do not see at this moment how
one can show the above uniqueness result by her algorithm.

Some of our statements  (after some 
identifications) can be compared with some results  of A. Pichon. 
E.g. our last formula from \ref{pr} can be compared with 
Proposition 3 \cite{Pi1}. 

\vspace{2mm}

\noindent Now we return back to our main theorem and its corollaries.

\vspace{2mm}

\noindent {\bf Corollary 1.} {\em Assume that $g(x,y,z)=f(x,y)+z^n$ is a 
suspension hypersurface singularity with $f$ irreducible, and not of type
described in the pathological cases S1 and S2 (cf. \ref{S1} and \ref{S2}).
Then the link $L_X$ of $(X,0)=(\{g=0\},0)$ 
determines completely the following data:

(1) The embedded topological type of $(X,0)$ (i.e. the embedding 
$L_X\subset S^5$), in particular, the Milnor fibration
and all the homological package derived from it. 

(2) All the equivariant Hodge numbers associated with the vanishing
 cohomology of $g$, in particular, the geometric genus of $(X,0)$. 

(3) The multiplicity of $g$.  

\vspace{1mm}

In particular, if $L_X$ is a rational homology sphere, then $L_X$ 
determines (1) (2) and (3). If $g$ is a pathological case listed in
S1 or S2, then $L_X$ together with the Milnor number of $g$ determines 
completely (1) (2) and (3).}

\vspace{2mm} 

Indeed, once we have the Newton pairs of $f$ and the integer $n$, then
the proof involves the description of the corresponding invariants
for plane curve singularities and  different ``Sebastiani-Thom type'' formulae,
 for more details  see e.g. \cite{AGV,SSS,NI,NII,NRev}. In fact, recently in 
\cite{NNIII}, the geometric genus (together with the Milnor number and the
signature of the Milnor fiber) was computed in terms of the 
Seiberg-Witten invariant of the link, provided that the link is a rational
homology sphere.  

It is well-known that the Milnor number of $g$ can be determined from
the embedded topological type of $g$ (a fact firstly noticed by Teissier,
see also \cite{YauConj}). Therefore, we get:

\vspace{2mm}

\noindent {\bf Corollary 2.} [Zariski's conjecture for this family]
 {\em The multiplicity  of $g=f(x,y)+z^n$ ($f$ irreducible) is determined by 
the embedded topological type of $g$. }

\vspace{2mm}

\noindent 
Corollary 1 (see also Theorem 2 below) motivates the following  conjecture.

\vspace{2mm}

\noindent {\bf Conjecture.} {\em Let $g:(\bc^3,0)\to 
(\bc,0)$  be an isolated  hypersurface singularity whose link $L_X$  is a 
rational homology sphere. Then the fundamental group of the link
characterizes completely the embedded topological type,
the equivariant Hodge numbers and the multiplicity of $g$. 
More generally, if the link of  a $\bfq$-Gorenstein
singularity $(X,0)$ is a rational homology sphere, then the  multiplicity 
of $(X,0)$ is determined by the (oriented) homeomorphism type of the link.}

\vspace{2mm}

Notice that the last general property is true for rational \cite{Artin,Artin2}
and elliptic  \cite{Lauferme,Yau,NInv} singularities.
The  first part can be verified in the following cases.
If $g=f+z^n$ with $f$ irreducible,  then  Conjecture is true
by Theorem 1 above (see also 
\cite{NePl} for the relation between $L_X$ and its fundamental group).
If $f$ is arbitrary but $n=2$, then it is true by \cite{Lauferntwo}.

For weighted homogeneous singularities 
the next theorem answers positively (since its proof is very short,
we decided to put the whole statement at the end of this introduction).

\vspace{2mm}

\noindent {\bf Theorem 2.} {\em The above conjecture is true for any 
weighted homogeneous hypersurface  singularity.}

\vspace{2mm}

\noindent Indeed, the Poincar\'e polynomial of the singularity can be
determined from the link by \cite{Dolg} or \cite{Pinkh}
(see also \cite{AbCov}).  Then, by a recent result of Ebeling \cite{Ebe}
(Theorem 1) follows that the characteristic polynomial of the  algebraic
monodromy can be recovered from the link $L_X$. (For this notice that
$\psi_A$ used by Ebeling is also link-invariant.) Then, by \cite{XY} 
(cf. also with \cite{OW}),  we get  that the weights, multiplicity
 and the embedded topological type is determined by $L_X$.
The statement about the Hodge data follows from \cite{St}. 

\section{Resolution graphs associated with analytic  germs}

\subsection{The embedded resolution.}\label{1.1} 
Let $(X,x)$ be a normal surface singularity and 
let $f:(X,x)\to (\bfc,0)$ be the germ of an analytic function.

An embedded resolution  $\phi:({\mathcal  Y},D)\to (U,f^{-1}(0))$ of
$(f^{-1}(0),x)\subset (X,x)$  is characterized
by the following properties.  There is a sufficiently small neighborhood 
$U$ of $x$ in $X$, smooth analytic manifold ${\mathcal  Y}$, and analytic
 proper map $\phi:{\mathcal  Y}\to U$ such that:

1)\ if $E=\phi^{-1}(x)$, then the restriction $\phi|_{{\mathcal  Y}
\setminus E}:
{\mathcal  Y}\setminus E\to U\setminus\{x\}$ is biholomorphic, 
and ${\mathcal  Y}\setminus E$ is dense in ${\mathcal  Y}$;

2)\ $D=\phi^{-1}(f^{-1}(0))$ 
is a divisor with only normal crossing singularities.

\noindent  $E$ is called the  exceptional divisor of $\phi$.
Let $\cup_{w\in \calw}E_w$ be its decomposition  into irreducible divisors.
The closure $S$ of $\phi^{-1}(f^{-1}(0)\setminus\{0\})$ is called the 
strict transform of $f^{-1}(0)$.
Let $\cup_{a\in \cala}S_a$ be its irreducible decomposition.
Obviously,  $D=E\cup S$. 

In the sequel  we will assume that $\calw\not=\varnothing$, 
any two irreducible components of $E$ have at most one intersection point, 
and  no irreducible exceptional divisor  has a self--intersection points.
This can always be realized by some additional blow ups.

\subsection{The embedded resolution graph.}\label{1.2}
We construct the embedded resolution graph  $\G(X,f)$ of the pair $(X,f)$,
associated with a fixed  resolution $\phi$,  as follows. Its vertices
$\calv=\calw\coprod\cala$
 consist of the nonarrowhead vertices $\calw$  corresponding to the
irreducible exceptional divisors, and arrowhead vertices $\cala$
correponding to the irreducible divisors of the strict transform $S$.
If two irreducible divisors corresponding to $v_1,v_2\in \calv$
have an intersection point then $(v_1,v_2)\ (=(v_2,v_1))$ is an edge of
$\G(X,f)$. 

The graph $\G(X,f)$ is decorated as follows. Any  $w\in \calw$ is 
decorated by  the self--intersection   $e_w:=E_w\cdot E_w$ and 
genus $g_w$ of $E_w$; and any $v\in \calv$ by 
the  multiplicity $m_v$ of $f$. More precisely,  for any $v\in\calv$,
 let $m_v$ be the vanishing order of $f\circ \phi$ along the
irreducible divisor corresponding to $v$. For example, if $f$ defines
an isolated singularity, then $m_a=1$ for any $a\in\cala$.

In all our graph--diagrams, we put the multiplicities in parentheses
(e.g.: (3)) and the genera in brackets (e.g.: [3]), with the convention 
that we  omit [0]. 

\subsection{The resolution of $(X,x)$.}\label{1.4} We 
say that $\phi:{\mathcal  Y}\to U$ is a resolution of $(X,x)$ if 
${\mathcal  Y}$ is a smooth analytic manifold, $U$ a neighbourhood of 
$x$ in $X$, 
$\phi$ is a proper analytic map, such that 
${\mathcal  Y}\setminus E$ (where $E=\phi^{-1}(x)$) is dense in
 ${\mathcal  Y}$ and
the restriction $\phi|_{{\mathcal  Y}\setminus E}:{\mathcal  Y}\setminus E\to
 U\setminus \{x\}$ is biholomorphic.

We codify the topology and the combinatorics of the pair $({\mathcal Y},E)$
in the dual resolution graph $\G(X)$ of $(X,x)$ associated with $\phi$. 
If the divisor $E$ is not a normal crossing divisor, then this codification 
can be slightly complicated, so in the sequel we will assume that the 
irreducible components of $E$ are smooth and intersect each other 
transversally, the irreducible components have no self--intersection points,
and there is no intersection point which is contained in 
more than two components.
In this case, similarly as in the situation of the embedded resolution,
the vertices of the dual graph
correspond to the irreducible components of $E$, the edges 
to the intersections of these components, and each vertex $w$ carries two
decorations: the genus $g_w$ of $E_w$ and the self--intersection
$E_w\cdot E_w$. 

Actually, one can obtain 
a possible graph $\G(X)$ from any $\G(X,f)$ by deleting 
all the arrows and multiplicities of the graph $\G(X,f)$.

An (embedded) resolution graph with the 
above restrictions is called {\em good} by some authors.

\subsection{Definitions.}\label{dn}   Let 
$\G$ be a decorated graph (as above, with or without arrowheads). 
For any $w\in \calw$, we denote by $\calv_w$ the set of vertices
$v\in \calv$ adjacent to $w$. Its cardinality $\#\calv_w$ is called the 
degree $\delta_w$ of $w$.
A vertex $w\in\calw$ is called {\em rupture} vertex if either \ $g_w>0$ \  or \
$\delta_w\geq 3$. 
(Notice that it can happen that $w$ is a rupture vertex in some $\G(X,f)$, 
but not in the graph obtained from $\G(X,f)$ after one deletes its arrows.) 
A vertex is called {\em leaf} vertex if $\delta_w=1$.

\subsection{Some  properties of resolution graphs.}\label{1.3}\ 

(1)\ Notice that the combinatorics of the 
graph and the self-intersections codify completely the intersection matrix 
$I:=(E_w\cdot E_v)_{(w,v)\in\calw\times
\calw}$ of the irreducible components of $E$. This matrix is 
negative definite, see e.g. \cite{Mumford} page 230; \cite{La} page 49, or 
\cite{GRa}.
We write $\det(\G):=\det(-I)>0$. By convention, the determinant 
of the empty graph is 1. 

(2) The graphs $\G(X,f)$ and $\G(X)$ are connected (see  \cite{La}, 
or Zariski's Main Theorem, e.g. in \cite{Hartshorne}).

(3)\ The next identities connect the  self--intersections 
and multiplicities: 
$$e_wm_w+\sum_{v\in\calv_w}m_v=0  \ \ \mbox{for any $w\in\calw$}.$$

Obviously, $m_v>0$ for any $v\in{\mathcal  V}$, 
hence the set of multiplicities determine the 
self--intersections completely. 
Similarly, since the intersection 
matrix $I$ is negative definite,  these  relations  determine the 
multiplicities $\{m_w\}_{w\in\calw}$ in terms
of the self--intersections   $\{e_w\}_{w\in\calw}$
and the multiplicities $\{m_a\}_{a\in\cala}$.

Using matrices, the above set of relations can be written as follows.
Fix a total  ordering of the set $\calw$. Let ${\bf m}_\calw$
be the column vector with $|\calw|$ entries $\{m_w\}_{w\in\calw}$.
Similarly, define the column vector ${\bf m}_\cala$   with $|\calw|$ entries 
whose  $w^{th}$ entry is $\sum_{a\in \cala\cap\calv_w}m_a$. Then 
\begin{equation*}
I\cdot {\bf m}_\calw+{\bf m}_\cala=0.
\end{equation*}
If $\G= \G(X)$ {\em is a tree}, then
the inverse matrix $I^{-1}$ can be computed in terms of 
determinants of some subgraphs as follows. Consider two vertices 
$w_1, \ w_2\in \calw$ and the shortest path which connects them. 
Let $\G_{w_1w_2}$ be the maximal
(in general non-connected) subgraph of $\G$ which has no vertices on this 
path. Then the $(w_1,w_2)$-entry of $I^{-1}$ is   given by 
\begin{equation*}
I^{-1}_{w_1w_2}=-\det(\G_{w_1w_2})/\det(\G).
\end{equation*}

(4) There are many (embedded) resolutions, hence many (embedded) resolution
graphs too. Nevertheless, they are all connected by quadratic
modifications (i.e. blow up and/or down of $-1$-curves with $g=0$).
Notice that by the above conventions,
we can blow down a $-1$-curve $E_w$ with $g_w=0$ if and only if
$\delta_w\leq 2$. If the (embedded) resolution 
graph has no rational  $-1$-curve $E_w$ with $\delta_w\leq 2$, then
we say that it 
is {\em minimal}. There is a unique minimal (embedded) resolution graph
denoted by $\G^{min}(X)$ (resp.  $\G^{min}(X,f)$). 
(Obviously, deleting the arrows from $\G^{min}(X,f)$ it can happen that
we obtain a non-minimal resolution graph of $(X,x)$.)

\subsection{The link of $(X,x)$.}\label{3.1}  Fix an embedding
$(X,x)\subset ({\bf C}^N,0)$ for some $N$. Then, for sufficiently small
$\epsilon>0$, $L_X:=\{z\in X: |z|=\epsilon\}$ is a connected 
oriented differentiable manifold 
(independent of the different choices). 
It is called {\em the link of  $(X,x)$}.

From topological point of view, $L_X$ characterizes completely $(X,x)$.
Moreover, in the presence of a resolution $\phi$,
for $U$ sufficiently small, the inclusion
$E\hookrightarrow {\mathcal Y}$ admits a strong deformation retract and
the restriction of $\phi$ 
identifies  $\partial {\mathcal  Y}$ with $L_X$. This shows that 
$L_X$ is  the {\em plumbed manifold}  associated with
$\G(X)$ (for details, see  \cite{NePl}), i.e. 
$\G(X)$ determines completely  $L_X$.
The converse is also true:  Neumann in \cite{NePl} proved that 
the (minimal)  resolution graph $\G(X)$ 
is determined  by the oriented homeomorphism type of $L_X$. 

\subsection{Fact. The homology of $L_X$.}\label{3.3} \cite{Sch,HNK,Mumford} \
{\em $H_1(L_X,\bfz)\approx \coker(I)\oplus \bfz^{2g+c_{\Gamma}},$
where $g:=\sum_{w\in \calw}g_w$ and $c_{\Gamma}$ denotes the number of 
independent cycles in $\Gamma=\Gamma(X)$ (e.g. $c_{\Gamma}=0$ if
and only if  $\Gamma$ is a tree). }

\vspace{2mm}

In particular, $L_X$ is an integral  (resp. rational) homology sphere 
if and only if  $g=c_{\Gamma}=0$ and $\det(\G)=1$ (resp. $g=c_{\Gamma}=0$).

\subsection{Example. Irreducible plane curve singularities.}\label{1.6}
(see e.g. \cite{BrKn} or \cite{EN}) \
If $(X,x)$ is smooth, hence $\approx (\bc^2,0)$,
then $f$ is called {\em plane singularity}. 
In this case, the graph $\G(\bc^2,f)$ is a tree, 
and $g_w=0$ for any $w\in\calw$. 

In this article we are mainly interested in {\em irreducible} plane curve 
singularities (i.e. when $|\cala|=1$). Their equisingular type  (and link also)
is completely characterized by the set of {\em Newton pairs}
$\{(p_k,q_k)\}_{k=1}^s$ (see e.g.  \cite{EN}, page 49).
Here  $(p_k,q_k)=1$,  $p_k\geq 2$, $q_k\geq 1$ and  $q_1>p_1$.

The {\em minimal} embedded resolution graph  can be reconstructed 
from the Newton pairs as follows
(see  e.g. in \cite{BrKn,EN} or \cite{NII}). 
First determine 
$u_i^l$ and $v_i^l$ ($u^0_i,v^0_i\geq 1$, and $u^l_i,v^l_i\geq 2$
for $l>0$) from the continued fractions:
$$\frac{p_i}{q_i}=u_i^0-{1\over\displaystyle u_i^1-
{\strut 1\over \displaystyle \ddots-
{\strut 1\over u_i^{t_i}}}};\hspace{1cm}
\frac{q_i}{p_i}=v_i^0-{1\over\displaystyle v_i^1-
{\strut 1\over\displaystyle \ddots-
{\strut 1\over v_i^{r_i}}}}.
$$
Then $\G^{min}(\bc^2,f)$ is given by 

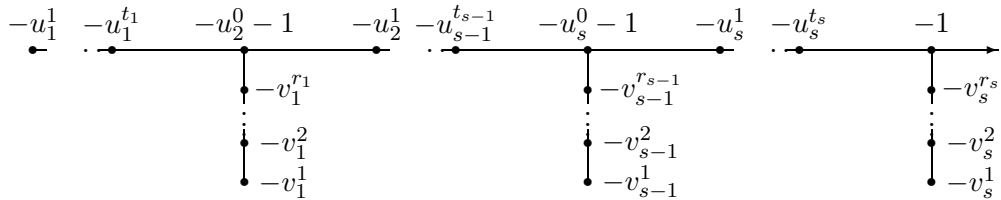
\vspace{5mm}

\begin{picture}(400,75)(20,0)
\put(20,65){\circle*{3}}
\put(50,65){\circle*{3}}
\put(100,65){\circle*{3}}
\put(150,65){\circle*{3}}
\put(180,65){\circle*{3}}
\put(230,65){\circle*{3}}
\put(280,65){\circle*{3}}
\put(310,65){\circle*{3}}
\put(360,65){\circle*{3}}
\put(100,50){\circle*{3}}
\put(100,30){\circle*{3}}
\put(100,15){\circle*{3}}
\put(230,50){\circle*{3}}
\put(230,30){\circle*{3}}
\put(230,15){\circle*{3}}
\put(360,50){\circle*{3}}
\put(360,30){\circle*{3}}
\put(360,15){\circle*{3}}
\put(20,65){\line(1,0){5}}
\put(45,65){\line(1,0){110}}
\put(175,65){\line(1,0){110}}
\put(305,65){\vector(1,0){80}}
\put(100,65){\line(0,-1){20}}
\put(100,35){\line(0,-1){20}}
\put(230,65){\line(0,-1){20}}
\put(230,35){\line(0,-1){20}}
\put(360,65){\line(0,-1){20}}
\put(360,35){\line(0,-1){20}}
\put(35,65){\makebox(0,0){\ldots}}
\put(165,65){\makebox(0,0){\ldots}}
\put(295,65){\makebox(0,0){\ldots}}
\put(100,40){\makebox(0,0){\vdots}}
\put(230,40){\makebox(0,0){\vdots}}
\put(360,40){\makebox(0,0){\vdots}}
\put(20,75){\makebox(0,0){$-u^1_1$}}
\put(50,75){\makebox(0,0){$-u^{t_1}_1$}}
\put(100,75){\makebox(0,0){$-u^0_2-1$}}
\put(150,75){\makebox(0,0){$-u^1_2$}}
\put(180,75){\makebox(0,0){$-u^{t_{s-1}}_{s-1}$}}
\put(230,75){\makebox(0,0){$-u^0_s-1$}}
\put(280,75){\makebox(0,0){$-u^1_s$}}
\put(310,75){\makebox(0,0){$-u^{t_s}_s$}}
\put(360,75){\makebox(0,0){$-1$}}
\put(115,50){\makebox(0,0){$-v^{r_1}_1$}}
\put(115,30){\makebox(0,0){$-v^2_1$}}
\put(115,15){\makebox(0,0){$-v^1_1$}}
\put(250,50){\makebox(0,0){$-v^{r_{s-1}}_{s-1}$}}
\put(250,30){\makebox(0,0){$-v^2_{s-1}$}}
\put(250,15){\makebox(0,0){$-v^1_{s-1}$}}
\put(375,50){\makebox(0,0){$-v^{r_s}_s$}}
\put(375,30){\makebox(0,0){$-v^2_s$}}
\put(375,15){\makebox(0,0){$-v^1_s$}}
\end{picture}

\noindent This has the following schematic form:

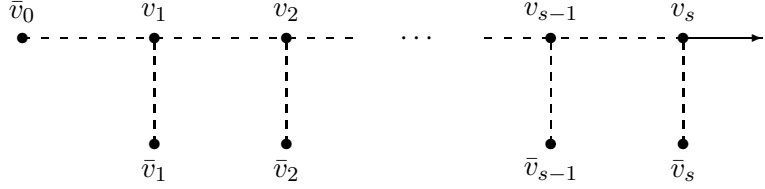
\begin{picture}(400,80)(0,0)
\put(50,60){\circle*{4}}
\put(100,60){\circle*{4}}
\put(150,60){\circle*{4}}
\put(250,60){\circle*{4}}
\put(300,60){\circle*{4}}
\put(100,20){\circle*{4}}
\put(150,20){\circle*{4}}
\put(250,20){\circle*{4}}
\put(300,20){\circle*{4}}
\put(50,70){\makebox(0,0){$\bar{v}_0$}}
\put(100,70){\makebox(0,0){$v_1$}}
\put(150,70){\makebox(0,0){$v_2$}}
\put(250,70){\makebox(0,0){$v_{s-1}$}}
\put(300,70){\makebox(0,0){$v_s$}}
\put(100,10){\makebox(0,0){$\bar{v}_1$}}
\put(150,10){\makebox(0,0){$\bar{v}_2$}}
\put(250,10){\makebox(0,0){$\bar{v}_{s-1}$}}
\put(300,10){\makebox(0,0){$\bar{v}_s$}}
\dashline{3}(50,60)(175,60)
\dashline{3}(225,60)(300,60)
\put(100,20){\dashbox{3}(0,40){}}
\put(150,20){\dashbox{3}(0,40){}}
\put(250,20){\dashbox{3}(0,40){}}
\put(300,20){\dashbox{3}(0,40){}}
\put(200,60){\makebox(0,0){$\cdots$}}
\put(300,60){\vector(1,0){30}}
\end{picture}

\noindent 
Here we emphasized only  those vertices $\{\bar{v}_k\}_{k=0}^s$
and $\{v_k\}_{k=1}^s$ which have degree  $\delta\not=2$.
We denote the  set of these vertices by ${\calw}^*$.
The dash-line between two such vertices replaces a string
\begin{picture}(100,10)(0,3)
\put(20,5){\circle*{4}}
\put(40,5){\circle*{4}}
\put(80,5){\circle*{4}}
\put(10,5){\line(1,0){40}}
\put(70,5){\line(1,0){20}}
\put(62,5){\makebox(0,0){$\cdots$}}
\end{picture}
. In our discussions below,
the corresponding self-intersections will be less important, but
the multiplicities of the vertices $v\in \calw^*$ will
be crucial. They can be easily described in terms 
of the integers $\{a_k\}_{k=1}^s$:
\begin{equation*}
a_1=q_1 \ \mbox{and } \  a_{k+1}=q_{k+1}+p_{k+1}p_ka_k\ \ 
\mbox{if \ $s-1\geq k\geq 1$}.
\tag{1}
\end{equation*}
Then again, $(p_k,a_k)=1$ for any $k$.
Clearly, the two sets  of pairs $\{(p_k,q_k)\}_{k=1}^s$ and
$\{(p_k,a_k)\}_{k=1}^s$  determine each other completely. 
In fact, the set of pairs 
$\{(p_k,a_k)\}_{k=1}^s$ constitutes the set of decoration of the so called
{\em  splice}, or {\em Eisenbud-Neumann}  diagram of $f$, cf. \cite{EN},
page 51. 
%
%
Then by \cite{EN}, section 10,  one has:
\begin{equation*}
\begin{array}{ll}
m_{v_k}=a_kp_kp_{k+1}\cdots p_s & \mbox{for $1\leq k\leq s$};\\
m_{\bar{v}_0}=p_1p_2\cdots p_s; & \\
m_{\bar{v}_k}=a_kp_{k+1}\cdots p_s & \mbox{for $1\leq k\leq s-1$};\\
m_{\bar{v}_s}=a_s.  &  \end{array}
\tag{2}
\end{equation*}

\subsection{ Example. Hirzebruch--Jung singularities.}\label{1.9}
(See \cite{Hir1,La,BPV}). 
For a normal surface singularity, the following conditions are equivalent.
If $(X,x)$ satisfies (one of) them, then it is called Hirzebruch--Jung 
singularity (where we prefer to include the smooth germ too).

(a)\ $(X,x)$ has a resolution graph $\G(X)$ which is a string with 
$g_w=0$ for any $w\in\calw$. (If the graph is minimal then 
it is either empty or  $e_w\leq -2$ for any $w$.)

(b)\ There is a finite proper map
$\rho:(X,x)\to ({\bf C}^2,0)$ such that reduced discriminant locus 
of $\rho$, in some local coordinates $(u,v)$ of $({\bf C}^2,0)$, is 
a subset of $\{uv=0\}$.

(c)\ $(X,x)$ is either smooth or is
isomorphic with exactly one of the ``model spaces'' 
$\{A_{n,q}\}_{n,q}$, where $A_{n,q}$ is the normalization of 
 $(\{xy^{n-q}+z^n=0\},0)$, where $0<q<n,\ (n,q)=1$.\\

For such $(X,x)$, $\pi_1(L_X)=\bz_n$, where 
$n=1$ if and only if $(X,x)$ is smooth (see e.g. \cite{BPV}), othervise
$n$ is the number from (c). 
Then, by \ref{3.3}, we also have $n=\det(\G(X))$ (fact which is valid for 
non-minimal graphs as well). 

In some of our applications we will need to recover this integer $n$ 
from the geometry of the map $\rho$ from (b). Consider the induced regular 
covering $\rho_{reg}:\rho^{-1}(\{uv\not=0\})\to \{uv\not=0\}$. Let 
$\pi_1=\bz^2$ be the fundamental group of $ \{uv\not=0\}$ generated by 
$e_1$ and $e_2$ representing two elementary loops around the axes $u$ and $v$.
Let $\rho_*:\pi_1\to G$ be the monodromy representation
of $\rho_{reg}$, and $\rho_*|\bz(e_i)$ be its  restriction to
$\bz(e_i)$, the subgroup generated by $e_i$ ($i=1,2$). 

\subsection{Lemma.}\label{n} \
$\bz_n\approx \ker(\rho_*)/\ker(\rho_*|\bz(e_1))\times \ker(\rho_*|\bz(e_2)).$
\begin{proof} This follows from the classification of the subgroups of
$\bz^2$, see e.g. \cite{BPV}, III.5. Another proof goes as follows.
It is clear that $\pi_1(\rho^{-1}(\{uv\not=0\}))=\ker(\rho_*)$,
and $\bz_n=\pi_1(L_X)$ is the quotient group of the previous group by the
subgroup generated by all the loops staying above $e_1$ and $e_2$. \end{proof}

\subsection{The ``model'' $X(a,b,N)$}\label{abN} From our point of view, 
it is more convenient to consider a bigger class of ``models'' instead
of $\{A_{n,q}\}_{n,q}$. More precisely, 
for any three strictly  positive integers $a,b$ and $N$, with $(a,b,N)=1$,
we define $(X,x)=
(X(a,b,N),x)$ as the unique singularity lying over the origin in 
the normalization of $(\{\alpha^a\beta^b+\gamma^N=0\},0)$. 
Let the germ $\gamma:(X(a,b,N),x)\to (\bfc,0)$ be 
induced by $(\alpha,\beta,\gamma)\mapsto \gamma$. 
In the sequel, we  give the embedded resolution
graph $\G(X,\gamma)$ of the germ $\gamma$ (for details, see \cite{NCC}). 

First,  consider  the unique $0\leq \lambda <N/(a,N)$  and 
$m_1\in{\bf N}$ with:
$$b+\lambda \cdot \frac{a}{(a,N)}=m_1\cdot \frac{N}{(a,N)}.$$
If $\lambda\not=0$, then consider the continued fraction:
$$\frac{N/(a,N)}{\lambda}=k_1-{1\over\displaystyle
k_2-{\strut 1\over\displaystyle\ddots
-{\strut 1\over k_s}}}, \ \ k_1,\ldots, k_s\geq 2.$$
Then the embedded resolution graph $\G(X,\gamma)$ of $\gamma$ is the 
following  string:

\begin{picture}(400,50)(-20,0)
\put(95,25){\makebox(0,0)[r]{$(\frac{a}{(a,N)})$}}
\put(40,25){\makebox(0,0)[r]{$St(a,b,N):$}}
\put(355,25){\makebox(0,0)[l]{$(\frac{b}{(b,N)})$}}
\put(150,35){\makebox(0,0){$-k_1$}}
\put(200,35){\makebox(0,0){$-k_2$}}
\put(300,35){\makebox(0,0){$-k_s$}}
\put(150,15){\makebox(0,0){$(m_1)$}}
\put(200,15){\makebox(0,0){$(m_2)$}}
\put(300,15){\makebox(0,0){$(m_s)$}}
\put(150,25){\circle*{4}}
\put(200,25){\circle*{4}}
\put(300,25){\circle*{4}}
\put(225,25){\vector(-1,0){120}}
\put(275,25){\vector(1,0){70}}
\put(250,25){\makebox(0,0){$\cdots$}}
\end{picture}

\noindent 
The arrow at the left (resp. right) hand side codifies the strict
transform of $\{\alpha=0\}$ (resp. of $\{\beta=0\}$). 
All vertices have genus $g_w=0$. The first vertex has
multiplicity $m_1$ given by the above congruence. 
Hence $m_2,\ldots, m_s$
can be easily computed using \ref{1.3}(3), namely:
$$-k_1m_1+\frac{a}{(a,N)}+m_2=0,\ \mbox{and}\
-k_im_i+m_{i-1}+m_{i+1}=0\ \mbox{for $i\geq 2$}.$$
This resolution resolves also the germs $\alpha$ and $\beta$
(induced by  the projection
$(\alpha,\beta,\gamma)\mapsto \alpha$, resp.
$(\alpha,\beta,\gamma)\mapsto \beta$). Indeed, as we already mentioned,
the strict transform of $\{\alpha=0\}$  
(resp. $\{\beta=0\}$) is  irreducible and it is exactly that strict transform 
component of $\gamma$ which is codified by the left (resp. right)
arrowhead of $St(a,b,N)$. The multiplicity of $\alpha $ (resp $\beta$)
along this strict transform component, or arrowhead, is $N/(a,N)$
(resp. $N/(b,N)$). Obviously, the multiplicity of $\alpha$ (resp. $\beta$)
on the right (resp. left) arrowhead is zero. 

Therefore, the embedded resolution graph $\G(X,\alpha^i\beta^j\gamma^k)$
($k>0$) of the germ  $\alpha^i\beta^j\gamma^k$ defined on $X$ 
can be deduced easily
from the above resolution graphs. It has the same shape, the same 
self--intersections and genera. Its multiplicity on the left arrowhead is
$(iN+ka)/(a,N)$, and on the right arrowhead $(jN+kb)/(b,N)$.

\vspace{2mm}

Finally, if $\lambda=0$, then the string $St(a,b,N)$
 has no vertices. In this case,  
$(X(a,b,N),x)$ is smooth and  the zero set of $\gamma$  (on $X$) has only 
a normal crossing singularity: in some local coordinates $(u,v)$ of $(X,x)$,
it can be represented as $\gamma=u^{a/(a,N)}v^{b/(b,N)}$.

\subsection{Remark.}\label{rm} \ 
Form the point of view of the above classification \ref{1.9}(c), 
$X(a,b,N)$ is an $A_{n,q}$--singularity,
where $n=N/(a,N)(b,N)$ and $q=\lambda/(b,N)$ 
(cf. e.g. \cite{BPV}, page 83-84). In  fact, $n$ can be deduced from 
\ref{n} as well. Indeed, in this case  $G=\bz_N$, $\rho_*$ is onto, hence 
$\ker(\rho_*)$ has index $N$ in $\bz^2$. On the other hand the index of
$\ker(\rho_*|\bz(e_i))$ in $\bz$ is $N/(a,N)$ for $i=1 $ and $N/(b,N)$ 
for $i=2$. 
Hence their product has index $N^2/(a,N)(b,N)$ in $\bz^2$. 
In particular, $\det(\G(X))$ is  $n=N/(a,N)(b,N)$. 

\subsection{Example.}\label{HJHJ} \ Consider 
the positive integers $q$,
$p$, $N$, $r$ and $P$ with $(q,p)=1$ and $(N,P)=1$. Define also $a:=rp+q$.
Let $(X_1,0):=(\{z^q=xy^p; \ w^N=(z^ry)^P\},0)\subset (\bc^4,0)$. 
Then $(X_1,0)$ is irreducible, its normalization is a Hirzebruch-Jung 
singularity with 
$$\det(-I)=n, \ \mbox{where} \ \  n=\frac{Nq}{(N,r)(N,a)}.$$
The proof has two steps. First we show that we can assume that $P=1$. 
Indeed,  consider the space $Y_1:=\{z^q=xy^p, \ w=t^P, \ t^N=z^ry\}$. 
Then $X_1$ and $Y_1$ are birationally equivalent (eliminate $t$ from the
equation of $Y_1$), hence their normalizations are the same. On the other hand
$Y_1$ is isomorphic with $X_1':=\{z^q=xy^p,\ t^N=z^ry\}$ (eliminate $w$).

Now we will apply \ref{n}. 
There is  natural finite map $(X_1,0)\to (\bc^2,0)$ induced by the 
$(x,y)$-projection which satisfies \ref{1.9}(b). The Galois group $G$ of the 
induced regular covering can be identified with 
$$G=\{(\xi,\eta)\in \bc^*\times \bc^*\  :\  \xi^q=1, \eta^N=\xi^r\}.$$
Since $\rho_*$ is onto, $\ker(\rho_*)$ has index $Nq$ in $\bz^2$ (cf. \ref{n}).
It is easy to see that $\rho(e_1)=(\exp(2\pi i/q),\exp(2\pi ir/(Nq)))$, hence 
$\ker(\rho_*|\bz(e_1))=k_1\bz$ for $k_1=Nq/(N,r)$. Similarly,
 $\rho(e_2)=(\exp(2\pi ip/q),\exp(2\pi ia/(Nq)))$, hence 
$\ker(\rho_*|\bz(e_1))=k_2\bz$ for $k_2=Nq/(N,a)$. Hence the claim  follows 
from \ref{n}.

\section{The resolution graph of cyclic coverings}

\noindent In this  section we review   some  properties of cyclic
coverings and their resolution graphs. The reader is invited to consult 
\cite{NCC} for more details. 

\subsection{Cyclic coverings}\label{6.1} \ 
Let $(X,x)$ be a normal surface singularity and 
$f: (X,x)\to (\bfc,0)$ the germ of an analytic function. For simplicity 
{\em we assume that $f$ defines an isolated singularity at $x$}. 
For any integer $n\geq 1$, take $b: (\bfc , 0)\to (\bfc ,0)$
given by $z\longmapsto z^{n}$, and let $X_{f,n}$ be the normalization
of the fiber product $\{ (x',z)\in (X\times\bfc , x\times 0) : f(x')=z^{n}\}$.
The second projection $(x',z)\in X\times\bfc\longmapsto z\in\bfc$ 
induces an analytic map 
$X_{f,n}\to\bfc$, still denoted by $z$.
The first projection $(x',z)\longmapsto x'$ gives 
rise to a ramified cyclic $n$-covering 
$pr  :X_{f,n}\to X$, branched along $\finv$. Since $\finv$ has an isolated
 singular point at $x$, one can verify that 
there is only one (singular) point of $\xfn$ lying above $x\in X$.

We regard  $\bfz _n$ as the group of $n^{th}$-roots of unity 
$\{ \xi_k=e^{2\pi ik/n};$ $\ 0\le k\le n-1\}$.
Then $(x',z)\longmapsto (x',\xi_k z)$ induces a 
$\bfz _n$--Galois action of $\xfn$ over $X$.

\subsection{The (embedded) resolution graph of cyclic coverings}\label{6.2} \ 
It is known (see e.g. \cite{NCC}) that in general the graphs $\G(\xfn,z)$ and
$\G(\xfn)$ {\em cannot} be determined from the embedded resolution graph 
$\G(X,f)$ of $f$ and the integer $n$. Nevertheless, if $L_X$ is a rational 
homology sphere, this is possible. The algorithm for a plane curve singularity
$f$ is given in \cite{NRev} (based on idea of the
construction used in the book of Laufer \cite{La}). 
For other particular cases, see \cite{Artal,O,OW,Lauferntwo}. 
The general case is discussed in \cite{NCC}. 
This algorithm  can also be compared with some results of 
E. Hironaka about global  cyclic coverings.

\subsection{Theorem. The graph $\G(\xfn,z)$}\label{6.3} \cite{NRev,NCC} \
Let $f:(X,x)\to(\bc,0)$ be as above  and assume that $L_X$ is a rational 
homology sphere. Fix the an embedded resolution graph $\G(X,f)$. 
Let $\{m_v\}_{v\in \calv}$ denote the corresponding multiplicities of $f$. 
For any $w\in \calw(\G(X,f))$, it is convenient to write 
$$M_w:=\mbox{gcd}(m_v\, |\, v\in \calv_w\cup\{w\}).$$
Similarly, if $e=(w_1,w_2)$ is an edge of $\G(X,f)$, then we write
$m_e:=(m_{w_1},m_{w_2})$.

Next, fix an integer $n$, and consider the cyclic covering $(\xfn,z)$
as above. 
In order to eliminate any confusion between the decorations  of the 
graph $\G(X,f)$ and $\G(\xfn,z)$, we denote the multiplicities
(respectively the genera)
of $\G(\xfn,z)$ by $\m_{v'} $ (respectively $\g_{v'}$
for any $v'\in \calv(\G(\xfn,z))$. 

Then a possible embedded resolution graph $\G(\xfn,z)$ can be 
constructed as follows. 

\vspace{2mm}

The graph $\G(\xfn,z)$ can be considered as a ``covering'' 
$q:\G(\xfn,z)\to \G(X,f)$.

\vspace{2mm}
 
(a)\ Above $w\in\calw(\G(X,f))$ there are $(M_w,n)$ vertices 
$v'\in q^{-1}(w)$ of $\G(\xfn,z)$,
each with multiplicity $\m_{v'}=m_w/(m_w,n)$ and genus $\g_{v'}$, where:
$$2-2\g_{v'}=\frac{(2-\delta_w)(m_w,n)+\sum_{v\in{\calv}_w}\ (m_w,m_v,n)}
{(M_w,n)}.$$
(In fact, by \ref{1.3}(3) , if 
$\delta_w\leq 2$ in $\G(X,f)$, then $\g_{v'}=0$.)
The vertices in $q^{-1}(w)$ can be indexed by the  group $\bz_{(M_w,n)}$.

\vspace{2mm}

(b)\ An edge $e=(w_1,w_2)$  of $\G(X,f)$ (where $w_1,w_2\in \calw(\G(X,f))$)

\begin{picture}(300,30)(0,0)
\put(150,10){\circle*{3}}
\put(200,10){\circle*{3}}
\put(150,10){\line(1,0){50}}
\put(150,20){\makebox(0,0){$(m_{w_1})$}}
\put(200,20){\makebox(0,0){$(m_{w_2})$}}
\end{picture}

\noindent
is covered by $(m_e,n)$ copies of strings in $\G(\xfn,z)$, each of type
(cf. \ref{abN}): 
$$St(m_{w_1}/(m_e,n),m_{w_2}/(m_e,n),n/(m_e,n)).$$
These strings can be indexed by the  group $\bz_{(m_e,n)}$.
The arrowheads of the strings are
identified with the vertices $q^{-1}(w_1)$, respectively $q^{-1}(w_2)$,
via the natural morphisms
$\bz_{(m_e,n)}\to \bz_{(M_{w_1},n)}$, respectively
$\bz_{(m_e,n)}\to \bz_{(M_{w_2},n)}$. 
(If this string is empty, i.e. if 
in \ref{abN} $\lambda=0$, then above the edge 
$e$ we insert $(m_e,n)$ edges by the same procedure.) 

\vspace{2mm}

(c)\ An arrowhead of $\G(X,f)$

\begin{picture}(300,30)(0,0)
\put(150,10){\circle*{3}}
\put(150,10){\vector(1,0){50}}
\put(150,20){\makebox(0,0){$(m_{w})$}}
\put(210,10){\makebox(0,0){$(1)$}}
\end{picture}

\noindent
is covered by one string of type $St(m_w,1,n)$ (cf. \ref{abN}), 
whose right arrowhead will
remain an arrowhead of $\G(\xfn,z)$ with multiplicity 1, and its left arrowhead
is identified with the unique vertex above $w$.
(Similarly as above, if this string is empty, then above this arrow 
 we insert a unique arrow supported by the unique vertex staying above $w$.) 

\vspace{2mm}

(d)\ In this way, we obtain all the vertices, edges and arrowheads
of $\G(\xfn,z)$, and all the multiplicities of  $z$.
Moreover, by the description of the
strings (cf. \ref{abN}), one has all the self--intersections
of the vertices which are situated on the new strings. The
self--intersections of the vertices $q^{-1}(w)\ (w\in{\calw}(\G(X,f)))$
can be computed using \ref{1.3}(3) (from the multiplicities of $z$). 

\vspace{2mm}

\noindent The isomorphism type of the above ``covering'' graph is independent 
of the choice a different identification (of the corresponding 
 sets with cyclic  groups), cf. \cite{NCC}. 
If we drop the arrowheads and multiplicities of $\G(\xfn,z)$, we obtain
$\G(\xfn)$. The graphs $\G(\xfn,z)$ and $\G(\xfn)$, in general, 
are not minimal.

\subsection{Definition.}\label{6.4} Let $(X,f)$ and $n$ be as in \ref{6.3}.
Assume that in the above algorithm, we start with the {\em minimal} 
(good) embedded
resolution graph $\G(X,f)$ of $(X,f)$. Then the output graph of the 
algorithm (without any modification by any blow up or down) will be called
the {\em canonical} embedded resolution graph of $(\xfn,z)$. 
In the sequel, we denote it by $\G^{can}(\xfn,z)$. 
The name is motivated by \cite{Lauferntwo}, where Laufer proved that 
the above algorithm for a plane curve singularity $f$ provides exactly
the canonical resolution  of $\xfn$ in the sense of Zariski, 
provided that $n=2$. 

\section{The resolution graph of  $\{f(x,y)+z^n=0\}$.}

In this section we make the algorithm \ref{6.3} very explicit in the case 
when $(X,x)=(\bc^2,0)$ and $f$ is an irreducible plane curve singularity.
Clearly, in this case, $\xfn$ can be identified with the hypersurface 
singularity $\{f(x,y)+z^n=0\}$ and $z$ with  the natural map induced by
the $z$-projection. We will  assume that $n\geq 2$.  

In the sequel we will use the following notations as well. 
Recall that $\{(p_k,q_k)\}_{k=1}^s$ denotes the set of
Newton pairs of $f$, and the integers $\{a_k\}_{k=1}^s$ are
 defined  in \ref{1.6}(1). Then we define

$\bullet$ \  $d_k:=(n,p_{k+1}p_{k+2}\cdots p_s)$ for
$0\leq k\leq s-1$, and $d_s:=1$; 

$\bullet$ \ $h_k:=d_{k-1}/d_k=(p_k, n/d_k)$ 
and $p_k':=p_k/h_k$ for $1\leq k\leq s$;

$\bullet$ \ $\tilde{h}_k:=(a_k,n/d_k)$ 
and $ a_k':=a_k/\tilde{h}_k$ for $1\leq k\leq s$.

\vspace{2mm}

We start our discussion with the computation of the numbers $M_w$
($w\in\calw^*$)  and
the decorations $\m_{v'}$ and $\g_{v'}$ from  \ref{6.3}(a) 
for any $v'\in q^{-1}(\calw^*)$. Notice that because of the 
Galois action, for a fixed $w\in 
\calw^*(\G^{min}(\bc^2,f))$, the integers $\m_{v'}$ and $\g_{v'}$ do not depend
on the choice of $v'\in q^{-1}(w)$, but only on $w$. Therefore, sometimes we 
prefer  to denote them simply by $\m_w$ and $\g_w$. 

\subsection{Lemma.}\label{l1} {\em Assume that $\G^{min}(\bc^2,f)$ is the  
minimal embedded 
 resolution graph of the irreducible plane curve singularity $f$.
Then the following facts hold:
$$\begin{array}{ll} 
M_{\bar{v}_k}=m_{\bar{v}_k} & (0\leq k\leq s) \ \ \ 
(\mbox{see \ref{1.6}(2) for  $m_{\bar{v}_k}$})\\
M_{v_k}=p_{k+1}\cdots p_s & (1\leq k\leq s-1)\\ 
M_{v_s}=1.  &  
\end{array}$$
For any $1\leq k\leq s$, fix two integers \ $i_k$ and $j_k$ with
$a_ki_k+p_kj_k=1$.   Then the multiplicities of the three vertices from the set
$\calv_{v_k}$, modulo $m_{v_k}=a_kp_kp_{k+1}\cdots p_s$, are:}
$$-i_ka_kp_{k+1}\cdots p_s; \ -j_kp_kp_{k+1}\cdots p_s; \ \mbox{and}
\   p_{k+1}\cdots p_s \ \ (k\leq s-1)$$
$$-i_sa_s; \ -j_sp_s; \ \ \mbox{and} \ \ \ 1 \ \ (k=s).$$
\begin{proof} The first identity follows from \ref{1.3}(3).
For the other statement 
see \cite{Ne2} or the proof of (3.2) in \cite{NI}. \end{proof}

\noindent Using this,  the graph $\G^{can}(\xfn,z)$ has the following data:

\subsection{Corollary.}\label{l2} (cf. also with \cite{NNII}.) {\em \ 
Let $q:\G^{can}(\xfn,z)\to \G^{min}(\bc^2,f)$
be the  ``graph projection'' considered in the algorithm \ref{6.3}. 
Then 

(a) For any $f$ and $n$, $ \G^{can}(\xfn,z)$ is a tree with 
$$\begin{array}{llll}
\#q^{-1}(v_s)=1 & \hspace{2cm} &
\#q^{-1}(v_k)=h_{k+1}\cdots h_s \ & (1\leq k\leq s-1)\\
\#q^{-1}(\bar{v}_s)=\tilde{h}_s & &
\#q^{-1}(\bar{v}_k)=\tilde{h}_kh_{k+1}\cdots h_s \ &  (1\leq k\leq s-1) \\
\#q^{-1}(\bar{v}_0)=h_1\cdots h_s. & & & 
\end{array}$$

(b) 

\vspace{-7mm}

$$\begin{array}{ll}
\m_{\bar{v}_0}=p_1'p_2'\cdots p_s' & \\
\m_{\bar{v}_k}=a_k'p_{k+1}'\cdots p_s' & (1\leq k\leq s-1)\\
\m_{\bar{v}_s}=a_s' &  \\
\m_{v_k}=a_k'p_k'p_{k+1}'\cdots p_s' & (1\leq k\leq s).\\
   \end{array}$$

(c) 

\vspace{-7mm}

$$\begin{array}{ll}
\g_{\bar{v}_k}=0 & (0\leq k\leq s)\\
\g_{v_k}=(h_k-1)(\tilde{h}_k-1)/2 & (1\leq k\leq s).
\end{array}$$
In particular, 
the link of $\xfn$ is a rational homology sphere if and only if 
$(h_k-1)(\tilde{h}_k-1)=0$ for any $1\leq k\leq s$ (cf. \ref{3.3}).}
\begin{proof}  Use \ref{6.3} and \ref{l1}.
For the last case in (c) notice that  $2-2\g_{v'}$ (with 
$\g_{v'}=\g_{v_k}$) equals
$-(a_kp_k,n/d_k)+(a_k,n/d_k)+(p_k,n/d_k)+1$, hence the identity
follows from the definition of $h_k$ and $\tilde{h}_k$ and $(a_k,p_k)=1$.
Notice that the fact that $\Gc$ is a tree follows also from Durfee's 
theorem \cite{DuF}, since the algebraic monodromy of $f$ has finite order
\cite{Le}.\end{proof}

\vspace{2mm}

\noindent By the above discussions,
the graph $\G^{can}(\xfn,z)$ has the following schematic form (where the
 dash-lines replace strings as above, and we omit the genera and the 
self-intersections):

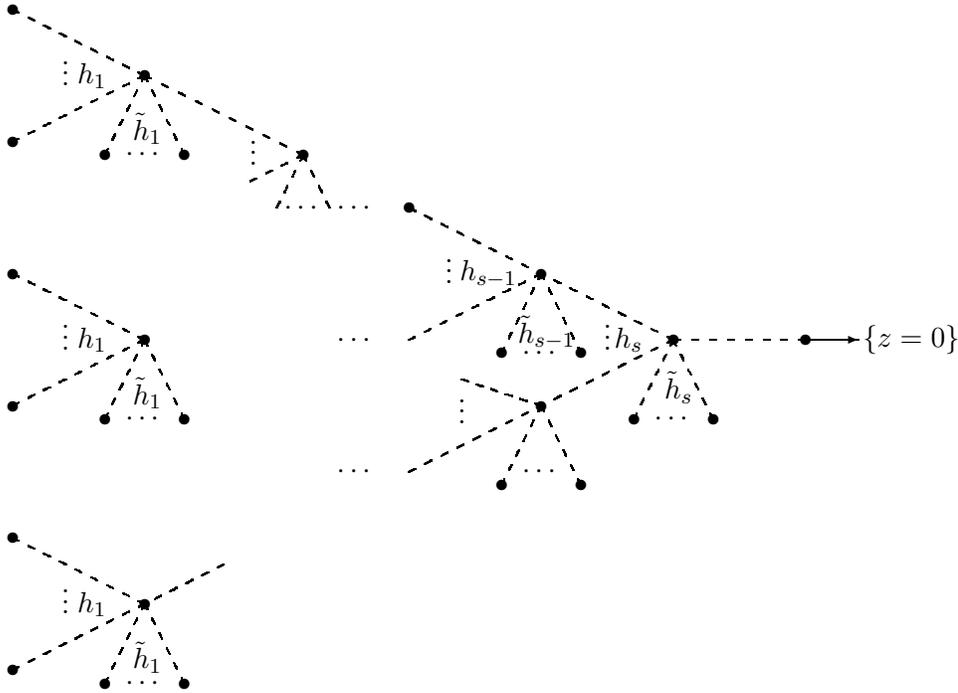
\begin{picture}(400,300)(30,-50)
\put(50,225){\circle*{4}}
\put(50,175){\circle*{4}}
\put(50,125){\circle*{4}}
\put(50,75){\circle*{4}}
\put(50,25){\circle*{4}}
\put(50,-25){\circle*{4}}
\put(100,200){\circle*{4}}
\put(100,100){\circle*{4}}
\put(100,0){\circle*{4}}
\dashline{3}(100,200)(50,225)
\dashline{3}(100,200)(50,175)
\dashline{3}(100,100)(50,125)
\dashline{3}(100,100)(50,75)
\dashline{3}(100,0)(50,25)
\dashline{3}(130,15)(50,-25)
\put(85,170){\circle*{4}}
\put(85,70){\circle*{4}}
\put(85,-30){\circle*{4}}
\put(115,170){\circle*{4}}
\put(115,70){\circle*{4}}
\put(115,-30){\circle*{4}}
\dashline{3}(100,200)(85,170)
\dashline{3}(100,100)(85,70)
\dashline{3}(100,0)(85,-30)
\dashline{3}(100,200)(115,170)
\dashline{3}(100,100)(115,70)
\dashline{3}(100,0)(115,-30)
\put(160,170){\circle*{4}}
\dashline{3}(160,170)(100,200)
\dashline{3}(160,170)(140,160)
\dashline{3}(160,170)(150,150)
\dashline{3}(160,170)(170,150)
\put(200,150){\circle*{4}}
\put(250,125){\circle*{4}}
\put(300,100){\circle*{4}}
\put(250,75){\circle*{4}}
\put(350,100){\circle*{4}}
\dashline{3}(200,150)(300,100)
\dashline{3}(200,50)(300,100)
\dashline{3}(250,125)(200,100)
\dashline{3}(250,125)(235,95)
\dashline{3}(250,125)(265,95)
\dashline{3}(250,75)(235,45)
\dashline{3}(250,75)(265,45)
\put(235,45){\circle*{4}}
\put(235,95){\circle*{4}}
\put(265,95){\circle*{4}}
\put(265,45){\circle*{4}}
\dashline{3}(300,100)(285,70)
\dashline{3}(300,100)(315,70)
\put(285,70){\circle*{4}}
\put(315,70){\circle*{4}}
\dashline{3}(300,100)(350,100)
\dashline{3}(250,75)(220,85)
\put(390,100){\makebox(0,0){$\{z=0\}$}}
\put(350,100){\vector(1,0){20}}
\put(70,204){\makebox(0,0){$\vdots$}}
\put(80,200){\makebox(0,0){$h_1$}}
\put(70,104){\makebox(0,0){$\vdots$}}
\put(80,100){\makebox(0,0){$h_1$}}
\put(70,4){\makebox(0,0){$\vdots$}}
\put(80,0){\makebox(0,0){$h_1$}}
\put(100,170){\makebox(0,0){$\cdots$}}
\put(101,180){\makebox(0,0){$\tilde{h}_1$}}
\put(100,70){\makebox(0,0){$\cdots$}}
\put(101,80){\makebox(0,0){$\tilde{h}_1$}}
\put(100,-30){\makebox(0,0){$\cdots$}}
\put(101,-20){\makebox(0,0){$\tilde{h}_1$}}
\put(141,174){\makebox(0,0){$\vdots$}}
\put(215,129){\makebox(0,0){$\vdots$}}
\put(230,125){\makebox(0,0){$h_{s-1}$}}
\put(220,76){\makebox(0,0){$\vdots$}}
\put(275,104){\makebox(0,0){$\vdots$}}
\put(283,100){\makebox(0,0){$h_{s}$}}
\put(250,50){\makebox(0,0){$\cdots$}}
\put(250,95){\makebox(0,0){$\cdots$}}
\put(252,103){\makebox(0,0){$\tilde{h}_{s-1}$}}
\put(300,70){\makebox(0,0){$\cdots$}}
\put(302,82){\makebox(0,0){$\tilde{h}_{s}$}}
\put(180,150){\makebox(0,0){$\cdots$}}
\put(160,150){\makebox(0,0){$\cdots$}}
\put(180,100){\makebox(0,0){$\cdots$}}
\put(180,50){\makebox(0,0){$\cdots$}}
\end{picture}

\subsection{Example.}\label{Seifert} Assume that $s=1$ and write $p=p_1$
and $a=a_1$. Take $n$ such that $h=(p,n)$ and $\tilde{h}=(a,n)=1$.
Then $\xfn$ can be identified with the Brieskorn hypersurface singularity
$\{(x,y,z)\in \bc^3:\ x^a+y^p+z^n=0\}$. Then the link is 
a Seifert 3-manifold with Seifert invariants: $a, \ a, \cdots, a, \ p/h, \ n/h$
($a$ appearing  $h$ times, hence all together there are
$h+2$ special fibers corresponding to the $h+2$ arms, cf. the above
graph-diagram).
These numbers also give (up to a sign) the determinants
of the corresponding arms of the  graph $\G(\xfn)$.
(For more details about Seifert manifolds and their plumbings, 
see e.g. \cite{JN,NeR}, or  \cite{NN}, section 6 for this special case; see
 also  \ref{S1} below). 

\subsection{The maximal strings of $\G^{can}(\xfn,z)$.}\label{str} \ 
The next goal is to compute the determinants of the maximal strings of
$\G^{can}(\xfn,z)$. For this, fix a vertex $w\in\calw^*(\G^{min}(\bc^2,f))$ 
and $v'\in q^{-1}(w)$. Consider the shortest path in $\G^{can}(\xfn,z)$ which
connects $v'$ and the arrowhead. 

If $w\not=v_s$, then on this path there is at least one rupture vertex  of
$\G^{can}(\xfn,z)$. Let $v''$ ($v''\not=v'$)
be the closest one to $v'$. If $w=v_k$ 
$(1\leq k\leq s-1)$, then let $\G(v')$ be the 
string which contains all the vertices between $v'$ and $v''$ (excluding
$v'$ and $v''$), and all the edges connecting them.
If $w=\bar{v}_k$  $(1\leq k\leq k)$, then $\G(v')$ is the string constructed
similarly, but at this time we include $v'$
and its  connecting  edge as well.
If $w=v_s$, then the above path is already a string.  Let $\G(v')$ be the
string which contains all the vertices between $v'$ and the arrowhead
 (excluding $v'$), and all the edges connecting them.

In this way we have  a codification of all the maximal strings of 
$\G^{can}(\xfn,z)$.
Notice also that the isomorphism type of the 
 string $\G(v')$ does not depend on the choice of 
$v'\in q^{-1}(w)$, but only on $w$. 
Therefore, sometimes it is preferable to denote this type  by $\G(w)$. 
Denote by $D(v')$ (or by  $D(w)$) the determinant $\det (\G(v'))$.
If $\G(w)=\varnothing$, then by definition $D(w)=1$ (cf. \ref{1.3}(1),
see also \ref{c1}).

\subsection{Proposition.}\label{pr} {\em \ Consider the vertices 
 $w\in \calw^* (\G^{min}(\bc^2,f))$ as above. 
Then $D(w)$ has the following values:}
$$\begin{array}{ll}
D(\bar{v}_0)=a_1' \ & \\
D(\bar{v}_k)=p_k' \ & (1\leq k\leq s)\\
D(v_s)=n/(h_s\tilde{h}_s) \ & \\
D(v_k)=n\cdot q_{k+1}/(d_{k-1}\tilde{h}_k\tilde{h}_{k+1})
 \ & (1\leq k\leq s-1)
\end{array}$$
\begin{proof} We start with the (``difficult'') case $D(v_k)$ 
($1\leq k\leq s-1$). Using the notations of \ref{1.6}, the maximal string in 
$\G^{min}(\bc^2,f)$ between  $v_k$ and $v_{k+1}$ has the following form:

\begin{picture}(400,50)(20,0)
\put(95,25){\makebox(0,0)[r]{($m_{v_k})$}}
\put(355,25){\makebox(0,0)[l]{$(m_{v_{k+1}})$}}
\put(150,35){\makebox(0,0){$-u^1_{k+1}$}}
\put(200,35){\makebox(0,0){$-u^2_{k+1}$}}
\put(300,35){\makebox(0,0){$-u^{t_{k+1}}_{k+1}$}}
\put(150,25){\circle*{4}}
\put(200,25){\circle*{4}}
\put(300,25){\circle*{4}}
\put(225,25){\vector(-1,0){120}}
\put(275,25){\vector(1,0){70}}
\put(250,25){\makebox(0,0){$\cdots$}}
\end{picture}

\noindent where $m_{v_k}=a_kp_k\cdots p_s$ and 
$m_{v_{k+1}}=a_{k+1}p_{k+1}\cdots p_s$, cf. \ref{1.6}(2). 
Moreover, $p_{k+1}/q_{k+1} =u^0_{k+1}-\lambda/q_{k+1}$, and the quotient
$q_{k+1}/\lambda$ gives the continued fraction $[u^1_{k+1},\cdots, u^{t_{k+1}}
_{k+1}]$. This can be identified with the string $St(a,b,N)$ in the
description of the model $X(a,b,N)$ in \ref{abN}. By this identification
$a=1$, $b=p_{k+1}$ and $N=q_{k+1}$. Therefore, the above string 
(without arrowheads) is the 
graph of the normalization of $\{z^{q_{k+1}}=xy^{p_{k+1}}\}$.
In this model, reading the 
multiplicities of the arrowheads associated with the coordinate functions, 
one gets for them $1$ and $p_{k+1}$ in the case $z$, and $0$ and $q_{k+1}$ 
in the case $y$. Therefore, $m_{v_k}$ and $m_{v_{k+1}}$ are the 
arrow-multiplicities of $(z^{p_ka_k}y)^{p_{k+1}\cdots p_s}$. 
In particular, the collection of graphs  $\{\G(v')\}_{v'\in q^{-1}(v_k)}$ 
is the (non-connected) graph of the normalization of 
$$X=\{(x,y,z,w)\,:\,
z^{q_{k+1}}=xy^{p_{k+1}}\ ;\ w^n=(z^{p_ka_k}y)^{p_{k+1}\cdots p_s}\}.$$
$X$ has $d_k=(n,p_{k+1}\cdots p_s)$ (isomorphic) irreducible components, 
number which agrees exactly with $\#q^{-1}(v_k)$. 
Hence, $D(v')$ is the graph of  the normalization of
$$X_1=\{z^{q_{k+1}}=xy^{p_{k+1}}\ ;\ w^{n/d_k}=(z^{p_ka_k}y)
^{(p_{k+1}\cdots p_s)/d_k}\}.$$
Then apply \ref{HJHJ} for $q=q_{k+1}$, $p=p_{k+1}$, $N=n/d_k$, $r=p_ka_k$,
$P=p_{k+1}\cdots p_s/d_k$ and $a=a_{k+1}$.

The other identities can be computed by a similar argument. But also notice 
that in all the other cases the corresponding maximal string contains a
leaf vertex of $\G^{can}(\xfn)$. Therefore, $D(v')$ can be 
identified with the corresponding Seifert invariant, similarly as in 
\ref{Seifert}. Then the first three identities also follow from \ref{Seifert}. 
\end{proof}

\subsection{Remark.}\label{c1} \ In $\Gc$ the following hold:

(1) If $w=\bar{v}_k$ ($0\leq k\leq s$) then $\G(w)\not = \varnothing$. Indeed,
$\G(w)$ contains at least as many vertices as the corresponding arm in 
$\G^{min}(\bc^2,f)$, which is clearly not empty. 

(2) The same argument is valid for any $\G(v_k)$  $(1\leq k\leq s-1)$ 
provided that $q_{k+1}>1$. In fact, for such $w=v_k$, $\G(v_k)=\varnothing$
if and only if $q_{k+1}=1$ and $n=d_{k-1}\tilde{h}_k\tilde{h}_{k+1}$. 

(3) $\G(v_s)=\varnothing$ if and only if $n=h_s\tilde{h}_s$.

\vspace{2mm}

Here appears a  natural question: is it  possible to distinguish the arms 
$\G(v')$ ($v'\in q^{-1}(\bar{v}_0)$) from the arms of type
$\G(v')$ ($v'\in q^{-1}(\bar{v}_1)$) ? The next corollary says that if 
$\g_{v_1}=0$ then already their  determinants are different:

\subsection{Corollary.}\label{c2} 

(a) If $D(\bar{v}_0)=D(\bar{v}_1)$ then $D(\bar{v}_0)=D(\bar{v}_1)=1$
and $\g_{v_1}\not=0$.

(b) If $D(\bar{v}_s)=D(v_s)$ then $D(\bar{v}_s)=D(v_s)=1$. 
\begin{proof} (a) If $D(\bar{v}_0)=D(\bar{v}_1)$, then 
$a_1/\tilde{h}_1=p_1/h_1$ by \ref{pr}. Since $(a_1,p_1)=1$, one  gets 
$a_1/\tilde{h}_1=p_1/h_1=1$. But then $h_1\geq 2$ and $\tilde{h}_1
\geq 2$ since $a_1=q_1>p_1\geq 2$. (b) Similarly by \ref{pr}
one has $p_s/h_s=n/h_s\tilde{h}_s$. But this two numbers are also 
relatively prime.
\end{proof}

\subsection{The subgraphs $\G_{\pm}(v_k)$.}\label{pm} \ Above we discussed 
the case of maximal strings  of $\G^{can}(\xfn,z)$. Obviously, one can 
consider the determinants of 
much bigger subgraphs delimited  by different rupture vertices. 
In this way one obtains a large number of rather subtle
invariants of this graph. Nevertheless,  
in order to recover the Newton pairs of $f$ and the integer $n$ from this 
graph, it is enough to consider only a restrictive sub-family of them.

Let us fix an integer $k$ $(1\leq k\leq s)$.
Consider the maximal subgraph of $\G^{can}(\xfn,z)$ which does not contain 
any  vertex from the set $q^{-1}(v_k)$. 
It has many connected components. 
The component which supports the arrowhead of $\G^{can}(\xfn,z)$ is denoted by
$\G_+(v_k)$. There are $\tilde{h}_kh_{k+1}\cdots h_s$ more 
components (isomorphic to each other), which contain
vertices above $\bar{v}_k$. They are strings of type $\G(\bar{v}_k)$
(cf. \ref{l2}(a) and  \ref{str}). Finally, there are $h_k\cdots h_s$ 
isomorphic components containing
vertices above $\bar{v}_0$. We denote such a component by $\G_-(v_k)$. 
$D_{\pm}(v_k)$ denotes $\det(\G_{\pm}(v_k))$. 

Obviously, $\G_-(v_1)=\G(\bar{v}_0)$, and $\G_+(v_s)=\G(v_s)$ whose 
determinants are computed in \ref{pr}. 

\subsection{Proposition.}\label{dd} {\em Assume that $s\geq 2$.

(a) $D_-(v_2)=(a_1')^{h_1-1}\cdot (p_1')^{\tilde{h}_1-1}\cdot a_2'$.

(b) $D_+(v_{s-1})=n\cdot D(v_{s-1})^{h_s-1}\cdot
D(\bar{v}_s)^{\tilde{h}_s-1}/(h_sh_{s-1}\tilde{h}_{s-1})$.}
\begin{proof} (a) Fix a vertex $v'\in q^{-1}(v_2)$ and one of the graphs
$\G_-(v')$. Let $w_1$ be its unique rupture point, and $w_2$ denote that 
vertex  which was connected by an edge  with $v'$ in $\G^{can}(\xfn,z)$.
(If $\G(v_1)=\varnothing$ then $w_1=w_2$, but the proof is valid in this
case as well.)
We put back on the vertices of  $\G_-(v')$ the multiplicities of 
$\G^{can}(\xfn,z)$. They will form a compatible set  (i.e. will satisfy
\ref{1.3}(3)) provided that we put on $w_2$ an arrow with multiplicity 
$\m_{v'}=\m_{v_2}$. 
This graph with arrowhead and multiplicities has the following
schematic form:

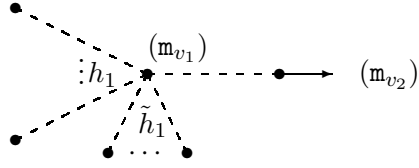
\begin{picture}(400,80)(130,55)
\put(300,100){\circle*{4}}
\put(250,125){\circle*{4}}
\put(250,75){\circle*{4}}
\put(350,100){\circle*{4}}
\dashline{3}(250,125)(300,100)
\dashline{3}(250,75)(300,100)
\dashline{3}(300,100)(285,70)
\dashline{3}(300,100)(315,70)
\put(285,70){\circle*{4}}
\put(315,70){\circle*{4}}
\dashline{3}(300,100)(350,100)
\put(380,100){\makebox(0,0)[l]{$(\m_{v_2})$}}
\put(300,110){\makebox(0,0)[l]{$(\m_{v_1})$}}
\put(350,100){\vector(1,0){20}}
\put(275,104){\makebox(0,0){$\vdots$}}
\put(283,100){\makebox(0,0){$h_{1}$}}
\put(300,70){\makebox(0,0){$\cdots$}}
\put(302,82){\makebox(0,0){$\tilde{h}_{1}$}}
\end{picture}

\noindent 
Notice that $\G_-(v')\setminus \{w_1\}$ has $h_1+\tilde{h}_1+1$ connected 
components, $h_1$ of type $\G(\bar{v}_0)$, $\tilde{h}_1$ of type 
$\G(\bar{v}_1)$, and one of type $\G(v_1)$. Therefore, 
by \ref{1.3}(3) one gets:
$$\frac{\m_{v_1}}{\m_{v_2}}=-I^{-1}_{w_1w_2}=\frac{D(\bar{v}_0)^{h_1}\cdot 
D(\bar{v}_1)^{\tilde{h}_1}}{D_-(v_2)}.$$
Now, use \ref{l2}(b) and \ref{pr}.  For part (b) we proceed similarly, 
but now with the graph $\G_+(v_{s-1})$. Its schematic form, together with the 
multiplicities of $\G^{can}(\xfn,z)$, is 

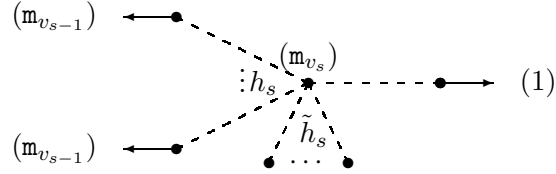
\begin{picture}(400,80)(130,55)
\put(250,125){\circle*{4}}
\put(250,75){\circle*{4}}
\put(350,100){\circle*{4}}
\dashline{3}(250,125)(300,100)
\dashline{3}(250,75)(300,100)
\put(275,104){\makebox(0,0){$\vdots$}}
\put(283,100){\makebox(0,0){$h_{s}$}}
\put(300,100){\circle*{4}}
\put(350,100){\circle*{4}}
\dashline{3}(300,100)(285,70)
\dashline{3}(300,100)(315,70)
\put(285,70){\circle*{4}}
\put(315,70){\circle*{4}}
\dashline{3}(300,100)(350,100)
\put(220,75){\makebox(0,0)[r]{$(\m_{v_{s-1}})$}}
\put(220,125){\makebox(0,0)[r]{$(\m_{v_{s-1}})$}}
\put(380,100){\makebox(0,0)[l]{$(1)$}}
\put(300,110){\makebox(0,0){$(\m_{v_s})$}}
\put(350,100){\vector(1,0){20}}
\put(250,75){\vector(-1,0){20}}
\put(250,125){\vector(-1,0){20}}
\put(300,70){\makebox(0,0){$\cdots$}}
\put(302,82){\makebox(0,0){$\tilde{h}_{s}$}}
\end{picture}

\noindent 
If from this graph we delete its rupture point (and the arrows and 
multiplicities) then  we get the following connected components: $h_s$ of type
$\G(v_{s-1})$, one of type $\G(v_s)$, and $\tilde{h}_s$ of type $\G(\bar{v}_s)
$. Therefore, from \ref{1.3}(3), similarly as above, one gets:
$$\m_{v_s}=\frac{ D(v_{s-1})^{h_s}\cdot D(\bar{v}_s)^{\tilde{h}_s}}
{D_+(v_{s-1})}+ h_s\cdot 
\frac{D(v_{s-1})^{h_s-1}\cdot  D(\bar{v}_s)^{\tilde{h}_s}\cdot  D(v_{s})}
 {D_+(v_{s-1})}\cdot 
\m_{v_{s-1}}.$$
Then use  again \ref{l2}(b) and \ref{pr} (and $a_s=q_s+p_sp_{s-1}a_{s-1}$).
\end{proof}

\subsection{Remark.}\label{gen} The above formulae and proofs  can be easily 
generalized for the other subgraphs as well. For example, 
for $\{D_-(v_k)\}_{k\geq 2}$ one can prove 
(by computing $\m_{v_{k-1}}/\m_{v_k}$ by the above method)
 the following inductive formula:
$$\frac{D_-(v_k)}{a_k'}=\Big[ \frac{D_-(v_{k-1})}{a_{k-1}'}\Big]^{h_{k-1}}
\cdot (a_{k-1}')^{h_{k-1}-1}\cdot (p_{k-1}')^{\tilde{h}_{k-1}-1}.$$

\vspace{3mm}

\section{From $\Gc$ to $\G^{min}(\xfn)$}

\noindent 
Let $\Gm$ be  the {\em minimal} embedded resolution graph of $(\xfn,z)$.
This can be obtained from $\Gc$ by a sequence of blow downs (and without any
blow up).

\subsection{Proposition.}\label{ru} {\em \
All the rupture vertices  of $\Gc$ survive in $\Gm$ as rupture vertices
 (i.e. they are not blown down in the minimalization procedure, and in
$\Gm$ they still live as  rupture vertices).}

\begin{proof}
From \ref{1.9} follows that a string of type $\G(\bar{v}_k)$ 
($0\leq k\leq s$) is completely collapsed in the minimalization procedure
 if and only if its determinant $D(\bar{v}_k)$ equals 1. 

First 
we verify that all the rupture vertices above $v_1$ will survive (as rupture
vertices). Let $v'$ be one of them considered in $\G^{can}$. 
It supports $h_1$ strings of type
$\G(\bar{v}_0)$, $\tilde{h}_1$ strings of type $\G(\bar{v}_1)$ and another 
edge, denoted by $e$. 
Recall that $D(\bar{v}_0)=a_1'$ and $D(\bar{v}_k)=p_k'$, cf. \ref{pr}.
By \ref{c2}, if both $D(\bar{v}_0)$ and $D(\bar{v}_1)$  equal one, then
$\g_{v'}\not=0$. Hence $v'$ will be a rupture vertex in $\Gm$.

If $D(\bar{v}_0)\not=1$
but $D(\bar{v}_1)=p_1/h_1=1$ then the strings of type 
$\G(\bar{v}_0)$ will survive. Their number is $h_1=p_1\geq 2$.
Symmetrically, if   $D(\bar{v}_1)\not=1$
but $D(\bar{v}_0)=a_1/\tilde{h}_1=1$ then $\tilde{h}_1=a_1\geq 2$
strings of type $\G(\bar{v}_1)$ will survive.
If both determinants are greater than one, then all the strings will survive
with total number  $h_1+\tilde{h}_1\geq 2$. 
Since the arrowhead survives, and $\Gm$ is connected, the edge $e$ will 
survive  as well. Hence $v'$ has degree at least three in $\Gm$. 

By induction, we assume that for a fixed $k$,  all the rupture vertices 
 above any $v_i\in \calw^*$ survive for any $i\leq k-1$.
 We show that this is the case for $v_k$ as well. For this, fix an
arbitrary $v'\in q^{-1}(v_k)$. 

First notice that by the inductive step, the $h_k$ subgraphs $\G_-(v')$ will
 survive (they cannot be completely contracted since they contain rupture 
points that survive). Similarly as above, since the arrowhead survives,
the edge connecting $v'$ with $\G_+(v_k)$ will also survive. 
If $D(\bar{v}_k)=1$, then $h_k=p_k\geq 2$. If
 $D(\bar{v}_k)\not=1$, then all the graphs $\G(\bar{v}_k)$ will survive.
Hence, in any case $\delta_{v'}\geq 3$ in $\Gm$. 
\end{proof}

Now, recall that  $\G^{min}(\xfn)$ denotes the minimal (good) resolution graph
of $(\xfn,0)$. It can be obtained from $\Gm$ by deleting its arrowhead (and
all the multiplicities) and blowing down successively all the $(-1)$-curves 
with genus zero and {\em new}
 degree $\leq 2$. In fact, there is exactly one case
when after deleting the arrowhead of $\Gm$ we do not obtain a minimal
graph, and this is described completely in the next proposition.
In the sequel we refer to this ``pathological'' situation as the ``P-case''. 

\subsection{Proposition.}\label{nonmin} {\em Assume that 
by deleting the arrowhead
of $\Gm$ we obtain a non-minimal graph. Then $\Gm$ has the following
schematic form with the two left branches isomorphic and with $e\leq -3$
(we omit the multiplicities).
The rational $(-1)$-curve  is the unique vertex $v'=q^{-1}(v_s)$
(which survives in $\Gm$, cf. \ref{ru}). 

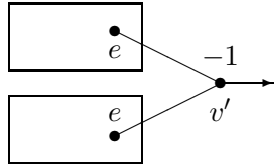
\begin{picture}(300,70)(0,15)
\put(200,50){\circle*{4}}
\put(160,70){\circle*{4}}
\put(160,30){\circle*{4}}
\put(200,50){\line(-2,1){40}}
\put(200,50){\line(-2,-1){40}}
\put(200,50){\vector(1,0){20}}
\put(120,20){\framebox(50,25){}}
\put(120,55){\framebox(50,25){}}
\put(200,60){\makebox(0,0){$-1$}}
\put(200,40){\makebox(0,0){$v'$}}
\put(160,62){\makebox(0,0){$e$}}
\put(160,38){\makebox(0,0){$e$}}
\end{picture}

This situation can happen if and only if $n=p_s=2$.

In this case, $\G^{min}(\xfn)$ is obtained from $\Gm$ by deleting its
arrowhead and blowing down $v'$. No other blow downs are necessary. 

Moreover, in this case, all the vertices of \, $\Gm$  have genus zero. }
\begin{proof}
If the graph obtained from $\Gm$ by deleting its arrowhead is not
minimal, then the vertex $v'$  in $\Gm$ which supports the arrowhead should be 
a $(-1)$ rational curve of degree 3 in $\Gm$. This can happen  only if this 
vertex  $v'$  is exactly the unique vertex $q^{-1}(v_s)$
(and $\g_{v_s}=0$). This also shows
that $\G(v_s)$ was collapsed in $\Gm$, 
hence $D(v_s)=n/h_s\tilde{h}_s=1$. Hence:
$$h_s\tilde{h}_s=n\geq 2\ \ \mbox{and} \ \ (h_s-1)(\tilde{h}_s-1)=\g_{v_s}=0.$$
Assume first that $h_s=1$ and $\tilde{h}_s>1$. Since $D(\bar{v}_s)=p_s>1$, the 
$\tilde{h}_s$ strings $\G(\bar{v}_s)$ are present in $\Gm$. This can happen if 
and  only if $\tilde{h}_s=2$ and $\G_-(v_s)$ is collapsed  completely in $\Gm$.
Since for $s\geq 2$ the rupture points $q^{-1}(v_1)$ survive in $\Gm$, 
this can 
happen  only if $s=1$ and $D(\bar{v}_0)=a_1/\tilde{h}_1=1$. This shows that 
$a_1=\tilde{h}_1=2$, which contradicts the inequality $a_1=q_1>p_1\geq 2$.

Therefore $\tilde{h}_s=1$ and $h_s>1$. Since the degree of $v'$ 
(in $\Gm$) is at least 
$1+h_s$ (hence $1+h_s\leq 3$), one gets $h_s=2$, and also the fact that 
the graphs of type $\G(\bar{v}_s)$ are collapsed in $\Gm$, hence $p_s/h_s=1$.
Therefore,  $n=h_s=p_s=2$. 

Then $\tilde{h}_k=1$, hence $\g_{v_k}=0$  for any $k$.

Finally, notice that $e\leq - 3$ (cf. the diagram) since after we blow down
$v'$ we get a subgraph of type \ \ 
\begin{picture}(60,10)(0,3)
\put(10,5){\circle*{4}}
\put(50,5){\circle*{4}}
\put(10,5){\line(1,0){40}}
\put(10,10){\makebox(0,0){$e+1$}}
\put(50,10){\makebox(0,0){$e+1$}}
\end{picture}
\ \ which must be  negative definite.

On the other hand, using the algorithm described in sections 3 and 4,
one can verify easily that if $n=p_s=2$ then  the above situation 
always occurs. \end{proof} 

\subsection{Remark.}\label{col} Assume that $\g_{v_k}=0$. Then  
 if a family of strings supported by any fixed $v'\in q^{-1}(v_k)$ 
is collapsed completely during the minimalization 
procedure, then the cardinality of this family (in spite of the fact that
it is missing in $\G^{min}(\xfn)$) can be determined, and it is one.
More precisely, if for $k\geq 1$, the $\tilde{h}_k$ graphs of type
$\G(\bar{v}_k)$ are completely collapsed, then $h_k=p_k\geq 2$. Then
by the genus-formula \ref{l2}(c) one gets  $\tilde{h}_k=1$.
Similarly, if $k=1$ and the $h_1$ graphs 
of type $\G(\bar{v}_0)$ are collapsed, then $h_1=1$. 

\vspace{2mm}

In order to recover the Newton pairs 
of $f$ and the integer $n$ from the graph $\G^{min}(\xfn)$, we need 
some information about some subgraphs 
 $\calg$ of $\G^{min}(\xfn)$ of the following type.  Each $\calg$ 
is a connected component of $\G^{min}(\xfn)\setminus \{v\}$ for 
some rupture vertex 
$v$ of $\G^{min}(\xfn)$, and  it contains exactly one rupture vertex  of 
$\G^{min}(\xfn)$. 
In a general setting their precise definition  is the following.
[We recall that the determinant of a (decorated) graph is the 
determinant of the negative of its intersection matrix (cf. \ref{1.3}).]

\subsection{Definitions.}\label{defs} Let $\G$ be a decorated tree
(with self-intersections and genera $\{\g_v\}_v$, 
without arrowheads and multiplicities). 
Assume that it has at least two rupture vertices. 

(1) Let $\G(\calr)$ be the minimal connected subgraph of $\G$ which contains
all the rupture vertices $\calr$  of $\G$.  Let $\call(\G(\calr))$ be the 
set of leaf vertices of $\G(\calr)$. For any $v\in \call(\G(\calr))$ let
$\calg(v)$ be the maximal connected subgraph of $\G$ which contains $v$ but 
contains  no other rupture vertex of $\G$. The determinant $\det(\calg(v))$
is denoted by $\cald(v)$. 

(2) For any $v\in \call(\G(\calr))$, let $v_{root}$ be the unique 
rupture vertex of $\G$ with the property that on the shortest  path
in $\G$ connecting $v$ and $v_{root}$ there are no other rupture vertices of
$\G$. Then clearly $v_{root} $ is adjacent with a certain  vertex of $\calg(v)$
(in fact $\calg(v)$ is one of the connected components of 
$\G\setminus \{v_{root}\}$). 

(3) For each rupture vertex $v\in \calr$, denote by $\calst(v)$
the set of maximal strings of $\G$ which are supported by $v$ (on one end)
and contain a leaf vertex of $\G$ (on the other end). More precisely, 
these strings are those connected components of $\G\setminus \G(\calr)$
which have an adjacent vertex with $v$ (in $\G$).
We write $\calst(v)$ as a disjoint union of its subsets $\{\calst_i(v)\}_{i\in
I(v)}$
which are the level sets of $\det:\calst(v)\to \bz$. We set $D_i:=\det(St)$
for $St\in \calst_i(v)$ and $\#_i:=\#\calst_i(v)$. Then we define 
$$\cald_{St}(v):=\left\{\begin{array}{ll}
\prod_{i\in I(v)} \ D_i^{\#_i}&\mbox{if $\calst(v)\not=\varnothing$}\\
1& \mbox{if $\calst(v)=\varnothing$,} \end{array}\right.
\ \  \ \ 
\cald_{St}^{red}(v):=\left\{\begin{array}{ll}
\prod_{i\in I(v)} \ D_i^{\#_i-1}&\mbox{if $\calst(v)\not=\varnothing$}\\
1& \mbox{if $\calst(v)=\varnothing$,} \end{array}\right.$$
and $\alpha(v)\in \bfq \cup\{\infty\}$ by 
$$\alpha(v)=\left\{\begin{array}{ll}
\prod_{i\in I(v)}\#_i & \ \mbox{if $\calst(v)\not=\varnothing$}\\
1 & \ \mbox{if $\calst(v)=\varnothing$ and $\g_v=0$}\\
\frac{2\g_v}{\delta_v-2}+1& \ \mbox{if $\calst(v)=\varnothing$ and $\g_v
\not=0$}.\end{array}\right.$$
[$\alpha(v)=\infty$ if and only if the degree $\delta_v$ of $v$ in $\G$ is 2,
$\calst(v)=\varnothing$ and $\g_v\not=0$.]

(4) For each $v\in \call(\G(\calr))$ we define the $\beta$-invariant by
$$\beta(v):=\frac{\cald(v)}{\cald_{St}(v)}\cdot
\frac{\alpha(v_{root})}{\alpha(v)}.$$

\subsection{}\label{nothing} In 
the next paragraphs we  apply these definitions for 
$\G=\G^{min}(\xfn)$. Here,  we prefer to regard $\G^{min}(\xfn)$
together with $\Gc$, as a minimalization of $\Gc$. In particular, we will
define subsets, subgraphs, etc.  in $\G^{min}(\xfn)$ 
as the images of well-defined
subsets, subgraphs,  etc. of $\Gc$  by the minimalization procedure.
(Of course,  in the next section will be a crucial task to recover some of 
these sets only from the abstract graph $\G^{min}(\xfn)$. 
The key result for this is the next \ref{iandii}.)

In order to avoid any confusion, for any subset of vertices of $\Gc$, 
we will denote by $\pi(A)$ the image of $A$ by the minimalization procedure.
Hence, $\pi(A)$ denotes those vertices of $\G^{min}(\xfn)$ which have
 ancestors in $A$, and survive in $\G^{min}(\xfn)$;
in some cases this set can be empty. 

\vspace{2mm}

The following facts follow easily from 
the structure-results proved in section 4 and the above propositions
\ref{ru} and \ref{nonmin}.

\subsection{Facts.}\label{lgr} {\em \ Assume that 
$\G=\G^{min}(\xfn)$ with $s\geq 2$. Then the following hold:

\vspace{2mm}

\noindent 
(a) The set $\call(\G(\calr))$ is the disjoint union of two sets
$\calr_1$ and $\call\calr_s$, where

(i) $\calr_1:=\pi(q^{-1}(v_1))$; 

(ii) $\call\calr_s:=\varnothing$ if $h_s>1$. Otherwise 
$\call\calr_s:=\pi(q^{-1}(v_s))$, 
the image by the minimalization procedure of the 
unique rupture vertex of $\Gc$ sitting above $v_s$.\\
(In both cases, by \ref{ru} and \ref{nonmin}, these sets are
subsets  of the rupture vertices  of $\G^{min}(\xfn)$.) 

\vspace{2mm}

\noindent 
(b) The subgraphs $\calg(v)$ for $v\in \call(\G(\calr))=\calr_1\cup\call
\calr_s$ 
(cf. part (a)) can be identified as follows:

(i) Assume that we are not in the ``P-case'' with $s=2$. 
For each $v\in q^{-1}(v_1)$
consider the unique subgraph of type $\G_-(v_{root})$  in  $\G^{can}(\xfn,z)$,
for some  $v_{root}\in q^{-1}(v_2)$,  which contains $v$. 
Then its image in $\G^{min}(\xfn)$ by the minimalization procedure is
  $\calg(v)$.

(ii) Assume that $h_s=1$. For $v= q^{-1}(v_s)$ 
consider $\G_+(v_{root})$ in  $\G^{can}(\xfn,z)$ with 
$v_{root}:=q^{-1}(v_{s-1})$.  Then its image in $\G^{min}(\xfn)$ by the 
minimalization procedure is $\calg(v)$.

\vspace{2mm}

\noindent  (c) $\alpha$ and $\beta$ are constant on $\calr_1$.}

\vspace{2mm}

\noindent 
(The motivation for the notation $\call\calr_s$ is the following: later
we will use 
the symbol $\calr_s$ for $\pi(q^{-1}(v_s))$; hence $\call\calr_s=\calr_s$ 
if $\pi(q^{-1}(v_s))$ is a ``leaf rupture vertex'', otherwise it is empty.)

\vspace{2mm}

The main point is that in (a), the cases (i) and (ii) can be distinguished by 
the genus  and  $\beta$-invariant. 

\subsection{Proposition.}\label{iandii} {\em Assume that $\G=\G^{min}(\xfn)$
with  $s\geq 2$.

(a) If  there exists at least one $v\in \call(\G(\calr))$ with $\g_v\not=0$, 
then $\calr_1=\{v\in \call(\G(\calr))\, :\, \g_v\not=0\}$ and 
$\call\calr_s=\{v \in \call(\G(\calr))\, :\, \g_v=0\}$
($\call\calr_s$ can be empty). 

(b) If  $\g_v=0$ for any $v\in\call(\G(\calr))$ and 
$\call\calr_s\not=\varnothing$,
 then $\beta(v)\in(0,\infty)$ and }
$$ \begin{array}{ll}
\beta(v)> 2 & \ \mbox{if $v\in \calr_1$},\\
\beta(v)\leq 1/2 & \ \mbox{if $v\in \call\calr_s$}.\end{array}$$
\begin{proof} (a) $\g_v$ is constant on $\calr_1$ and 
$\g_v=0$ for (the unique) $v\in\call \calr_s$ provided that 
$\call\calr_s\not=\varnothing
$, since in this case $h_s=1$ (cf. \ref{l2}). 

Now we prove (b). Since $h_s=1$ we can exclude the ``P-case''. 

First assume that $s\geq 3$. 

If $v\in \calr_1$ then by \ref{c2} and \ref{col} one has $\alpha(v)=h_1
\tilde{h}_1$. For $v_{root}$, analyzing the three different cases from the 
definition of $\alpha(v_{root})$, and using \ref{col} and the genus formula,
we get $\alpha(v_{root})=\tilde{h}_2$. On the other hand,  
$\cald(v)=(a_1')^{h_1-1}(p_1')^{\tilde{h}_1-1}a_2'$ (cf. \ref{dd}(a))
and $\cald_{St}(v)=(a_1')^{h_1}(p_1')^{\tilde{h}_1} $ (use \ref{pr} and
notice that if a string
is collapsed completely then its determinant is one).
Therefore, $\beta(v)=a_2/(a_1p_1)>p_2\geq 2$,  cf. \ref{1.6}(1).

If $v=\pi(q^{-1}(v_s))$ then $\alpha(v)=\tilde{h}_s$ (use $\g_s=0$, \ref{col}
and \ref{c2}). By similar argument as above,  $\alpha(v_{root})=
\tilde{h}_{s-1}$. By \ref{dd}(b) and $h_s=1$ one has  $\cald(v)=D_+(v_{s-1})=
n(p_s')^{\tilde{h}_s-1}/(h_{s-1}\tilde{h}_{s-1})$.  By \ref{pr},
$\cald_{St}(v)=(p_s')^{\tilde{h}_s}\cdot n/(h_s\tilde{h}_s)$. 
Therefore, using again $h_s=1$ one gets $\beta(v)=1/(h_{s-1}p_s)\leq 1/2$.

Assume that $s=2$ and let $v'_i=q^{-1}(v_i)$ ($i=1,2$). Then 
$\alpha(v_i')=h_i\tilde{h}_i$ ($i=1,2$). Hence the computation of $\beta(v_1')$
 is the same as above, and it gives $a_2/(a_1p_1)>2$. For $v_2'$ we have an 
additional $h_1$ and we get $\beta(v_2')=1/p_2\leq 1/2.$
\end{proof}

\section{From $\G^{min}(\xfn)$ back to $f$ and $n$}

\noindent Our final goal is to recover the Newton pairs of $f$ and the 
integer $n$ from the graph $\G^{min}(\xfn)$. In general, this is not 
possible. Nevertheless, by our main theorem, there are only two cases when
 such an  ambiguity appears. They are presented in the next subsections.

\subsection{Example. The S1-coincidence.}\label{S1} Assume that $(X,0)=
(x^{q_1}+y^{p_1}+z^n=0,0)$ is a Brieskorn singularity with $(q_1,p_1)=1$. 
Let us first analyze how one can recover the set of integers $\{q_1,p_1,n\}$
from the minimal resolution graph $\G$ of $(X,0)$. In this case,
the computation of the graph $\G$ from the integers $\{q_1,p_1,n\}$  is a 
classical, well-known fact (cf. also with our algorithm). 
The graph is either a string (with all genera zero) or a star-shaped 
graph (where only the central vertex might have a non-zero genus).
If $\G$ is a string, then $(X,0)$ is a Hirzebruch-Jung {\em hypersurface}
singularity. But there is only one family of such singularities, namely
the  $A_{q_1-1}$-singularities provided by the integers $\{q_1,2,2\}$.
In this case, $q_1$ is just the determinant of the string.

There is a rich literature of star-shaped graphs and Seifert 3-manifolds, and
also of their subclass given by Brieskorn hypersurface singularities. 
The reader is invited to consult \cite{OW}, section 3, case (I) (cf. also with 
\cite{JN,NeR}).

If one wants to recover the integers $\{q_1,p_1,n\}$, then one considers the 
set of strings $\calst(v)$ of the central vertex $v$. 
Recall \ref{defs}(3) for the  notations.  Then $\#I(v)\leq 3$. 
If $\#I(v)=3$ then  $\{q_1,p_1,n\}=\{D_1\#_2\#_3,D_2\#_1\#_3,D_3\#_1\#_2\}$.
If one $\calst_{i_0}$ is missing (empty) then $D_{i_0}=1$ and 
$\#_{i_0}$ can  determined from the genus of the central vertex 
(see e.g. our genus formula \ref{l2}(c) or \cite{OW} (3.5)). Hence the previous
procedure  still works.

Similarly, in our special situation $(q_1,p_1)=1$, one can show
that if two subsets $\calst_i$ are empty, then one can still recover
$\{q_1,p_1,n\}$ excepting only one case, namely when $\calst(v)$ consists
of only one string. In our terminology, this can happen  only when 
the string which supports the arrowhead survives and all the others are 
contracted (i.e. $p_1'=a_1'=1$). 
Similar ambiguity appears when $\calst(v)=\varnothing$. 

But all these ambiguity cases can be classified very precisely. Consider 
an identity of type $(h_1-1)(\tilde{h}_1-1)=2\g>0$ and an arbitrary positive 
integer $l$. Then the triplet $\{q_1,p_1,n\}=
\{\tilde{h}_1,h_1,\tilde{h}_1h_1l\}$ provides the following graph $\G$
(with $l-1$ $(-2)$-vertices):

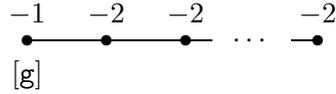
\begin{picture}(100,40)(-100,-10)
\put(20,10){\circle*{4}}
\put(50,10){\circle*{4}}
\put(80,10){\circle*{4}}
\put(130,10){\circle*{4}}
\put(20,10){\line(1,0){70}}
\put(130,10){\line(-1,0){10}}
\put(105,10){\makebox(0,0){$\cdots$}}
\put(20,-3){\makebox(0,0){$[\g]$}}
\put(20,20){\makebox(0,0){$-1$}}
\put(50,20){\makebox(0,0){$-2$}}
\put(80,20){\makebox(0,0){$-2$}}
\put(130,20){\makebox(0,0){$-2$}}
\end{picture}

\noindent Now fix $l>0$ and $\g>0$. Then, different triplets 
$\{\tilde{h}_1,h_1,\tilde{h}_1h_1l\}$ with $(h_1-1)(\tilde{h}_1-1)=2\g$,
$h_1>1$ and $\tilde{h}_1>1$ provide the same graph. 
[For example, $(3,7,21)$ and $(4,5,20)$ provide the same graph consisting 
of a vertex with $\g=6$ and self-intersection $-1$. ]

This is the only coincidence in the case of Brieskorn singularities with 
$(q_1,p_1)=1$. Obviously, this cannot happen if $\g=0$. 

\vspace{2mm}

\noindent {\em Addendum. Relation with the Milnor number.} Notice that in 
those cases when $\G$ fails to determine the integers $\{q_1,p_1,n\}$, 
$\G$ together with the Milnor number $\mu$ of the Brieskorn singularity
do determine $\{q_1,p_1,n\}$. Indeed, in the ``ambiguity cases'' one has
$(q_1,p_1,n)=(\tilde{h}_1,h_1,\tilde{h}_1h_1l)$, where
$2\g=(h_1-1)(\tilde{h}_1-1)>0$  and $l$ are 
readable from the graph. But $\mu=(h_1-1)(\tilde{h}_1-1)(\tilde{h}_1
h_1l-1)=2\g(\tilde{h}_1 h_1l-1)$. This determines $\tilde{h}_1h_1$, and
finally $h_1$ and $\tilde{h}_1$ (using the genus formula).

\subsection{Remark. The ``$z$-axis ambiguity''.}\label{zamb} Recall that by 
our general aim, we have to 
recover  the Newton pairs of $f$ and the integer $n$. In the Brieskorn 
case, after we recover the set 
$(q_1,p_1,n)$ we have to make a choice for the $z$-axis.
Recall that $(p_1,q_1)=1$. If 
$k$ integers among of $(p_1,q_1), \ (p_1,n)$ and $(q_1,n)$ equal 1,
then there are $k$ possibilities for the choice of the $z$-axis. 

\subsection{Example. The S2-coincidence.}\label{S2} The next coincidence 
appears when 
\begin{equation*}
\mbox{$s=2$, and $a_1'=p_1'=\tilde{h}_2=1$ (or equivalently,
$s=2$, $(n,a_2)=1$  and $(n,p_2)a_1p_1|n$\ ). }
\tag{$*$}
\end{equation*}
In this case clearly $q_1=\tilde{h}_1$ and $p_1=h_1$, but the $\tilde{h}_1$
strings of type $\G(\bar{v}_0)$ and the $h_1$ strings of type
 $\G(\bar{v}_1)$ are not visible 
on the  minimal graph  since their determinants are one, hence they are 
contracted. The graph $\G^{min}(\xfn) $ has the following schematic  form,
where  $\g_{v_1}>0$  and we omit the self-intersections:

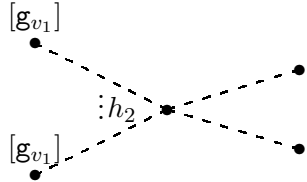
\begin{picture}(400,90)(130,55)
\put(300,100){\circle*{4}}
\put(250,125){\circle*{4}}
\put(250,75){\circle*{4}}
\put(350,115){\circle*{4}}
\put(350,85){\circle*{4}}
\dashline{3}(250,125)(300,100)
\dashline{3}(250,75)(300,100)
\dashline{3}(300,100)(350,115)
\dashline{3}(300,100)(350,85)
\put(250,85){\makebox(0,0){$[\g_{v_1}]$}}
\put(250,135){\makebox(0,0){$[\g_{v_1}]$}}
\put(275,104){\makebox(0,0){$\vdots$}}
\put(283,100){\makebox(0,0){$h_{2}$}}
\end{picture}

\noindent The strings that appear on the right correspond to $\G(v_2)$ and 
$\G(\bar{v}_2)$, but in general, we cannot decide which one is which.
From the graph we can read $h_2$ and the genus $\g_{v_1}
=(h_1-1)(\tilde{h}_1-1)/2$,
and of course, a lot of determinants.

Using $h_2, \ \tilde{h}_2$, and 
$D(v_1),\ D_-(v_2),\ D_+(v_1)$ and the set $\{D(v_2),D(\bar{v}_2)\}$,
we can also recover $a_2,\ q_2,  \ n/(h_1\tilde{h}_1),\ h_1\tilde{h}_1p_2$ and
the set $\{p_2,n\}$, where we cannot distinguish $p_2$ from $n$.

Notice that once we know $h_1\tilde{h}_1$, then using the genus formula
and $\tilde{h}_1=q_1>p_1=h_1$, we obtain $h_1$ and $\tilde{h}_1$ without
any ambiguity, hence (by the above equations) all the data. But for the three
``variables'' $h_1\tilde{h}_1,\ p_2,\ n$  we know only 
the values $n/(h_1\tilde{h}_1),
\ h_1\tilde{h}_1p_2$ and the set $\{p_2,n\}$. This, in general, has two 
possible solutions (which correspond by a permutation of  $p_2$ and $n$). 
If this is the case, 
then it might happen that there are two different realizations of the 
same graph for two different pairs $(f,n)$. But for this,  both solutions 
should provide positive integers as candidates for 
the Newton pairs and $n$. If this is not happening then the graph is 
uniquely realized (see Example 3 below). 

The complete discussion of all the cases when the above equations which 
involve  
$D(v_1)$, $D_-(v_2)$, $D_+(v_1)$ and the set $\{D(v_2),D(\bar{v}_2)\}$ 
associated with the  graph  provide exactly two ``good'' solutions for 
$(f,n)$ is long and tedious, so we decided  not to give it  here
(nevertheless we think that Example 3 illuminates completely the problem).
What is important is the fact that any graph 
(in this family) can be realized by at most 
two possible pairs $(f,n)$, and this coincidence in some cases really occurs.
(Moreover, given a pair $(f,n)$, or the graph of $\xfn$, one can write down 
easily the  possible candidate for the numerical data of 
$(f',n')$, the possible pair of $(f,n)$,  with the same graph.)

In the next examples we will write $\{(p_1,q_1),(p_2,q_2);n\}$ for the Newton
pairs of $f$ and the integer $n$. Recall that  ($*$) implies $p_1=h_1$ and
$q_1=\tilde{h}_1$. 

\vspace{2mm}

\noindent {\bf Example 1.} The two different solutions $(3,7)$ and $(4,5)$ for
the genus formula $(h_1-1)(\tilde{h}_1-1)=2\cdot 6$ can be completed 
to the following two sets of invariants: $\{(3,7),(20,1);21\}$ and
$\{(4,5),(21,1);20\}$. For them  the corresponding two graphs are the same:

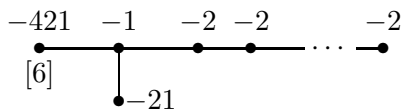
\begin{picture}(400,50)(0,20)
\put(100,50){\circle*{4}}
\put(130,50){\circle*{4}}
\put(160,50){\circle*{4}}
\put(180,50){\circle*{4}}
\put(230,50){\circle*{4}}
\put(130,30){\circle*{4}}
\put(100,50){\line(1,0){100}}
\put(230,50){\line(-1,0){10}}
\put(130,50){\line(0,-1){20}}
\put(100,40){\makebox(0,0){$[6]$}}
\put(100,60){\makebox(0,0){$-421$}}
\put(130,60){\makebox(0,0){$-1$}}
\put(160,60){\makebox(0,0){$-2$}}
\put(180,60){\makebox(0,0){$-2$}}
\put(230,60){\makebox(0,0){$-2$}}
\put(142,30){\makebox(0,0){$-21$}}
\put(210,50){\makebox(0,0){$\cdots$}}
\end{picture}

\noindent where the number of $(-2)$-curves is 19. Here $h_2=1$. 

\vspace{2mm}

\noindent {\bf Example 2.} If one wants examples with arbitrary $h_2$,
then one of the possibilities is the following:
one multiplies in the above data (of Example 1)  $p_2$ and $n$ by the wanted
$h_2$. E.g.  the data  $\{(3,7),(40,1);42\}$ and
$\{(4,5),(42,1);40\}$ provide the same $h_2=2$ and the same graph:

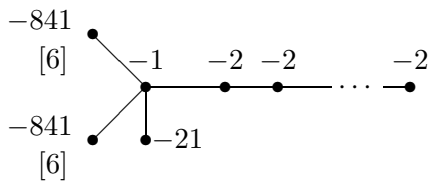
\begin{picture}(400,65)(0,15)
\put(110,30){\circle*{4}}
\put(110,70){\circle*{4}}
\put(130,50){\circle*{4}}
\put(160,50){\circle*{4}}
\put(180,50){\circle*{4}}
\put(230,50){\circle*{4}}
\put(130,30){\circle*{4}}
\put(130,50){\line(1,0){70}}
\put(110,70){\line(1,-1){20}}
\put(110,30){\line(1,1){20}}
\put(230,50){\line(-1,0){10}}
\put(130,50){\line(0,-1){20}}
\put(95,60){\makebox(0,0){$[6]$}}
\put(90,75){\makebox(0,0){$-841$}}
\put(95,20){\makebox(0,0){$[6]$}}
\put(90,35){\makebox(0,0){$-841$}}
\put(130,60){\makebox(0,0){$-1$}}
\put(160,60){\makebox(0,0){$-2$}}
\put(180,60){\makebox(0,0){$-2$}}
\put(230,60){\makebox(0,0){$-2$}}
\put(142,30){\makebox(0,0){$-21$}}
\put(210,50){\makebox(0,0){$\cdots$}}
\end{picture}

\noindent  where again, the number of $(-2)$-curves is 19. 

\vspace{2mm}

\noindent {\bf Example 3.} Assume the data  $\{(p_1,q_1),(p_2,q_2);n\}$
of $(f,n)$ satisfies ($*$), hence $p_1=h_1$ and $q_1=\tilde{h}_1$. 
If $(f,n)$ has ``a pair'' $(f_2,n_2)$ (with the same graph) then the 
data of $(f_2,n_2)$ has the form (cf. the above discussion)
$\{(x,y),(n,q_2),p_2\}$, where  $x$ and $y$ can be 
determined by the equations: $xy/p_2=p_1q_1/n$ and $(x-1)(y-1)=(p_1-1)(q_1-1)$.
It is easy to write down cases when this has no integral solutions.

E.g., the data $\{(2,3),(5,1);6\}$ satisfies ($*$), but it has no ``pair''.
Its minimal resolution graph can be realized
in a unique way in the form $f+z^n$ ($f$ irreducible) (cf. \ref{main}).

\vspace{2mm}

\noindent {\em Addendum. Relation with the Milnor number.} Even if the same
graph is realized by two different pairs $(f_1,n_1)$ and $(f_2,n_2)$, the 
corresponding Milnor numbers $\mu_i$ associated with the hypersurface
singularities $f_i+z^{n_i}$ ($i=1,2$) distinguish the two cases. 
This follows from the formula $\mu=[2\g_{v_1}p_2+(p_2-1)(a_2-1)](n-1)$.
Since $\g_{v_1}>0,\ a_2, np_2$ and $n+p_2$ are readable from the graph, 
this relation determines $p_2$, hence all the numerical data. 

\vspace{3mm}

\noindent 
Now we are ready to formulate and prove the main result of the article.

\subsection{Theorem.}\label{main} {\em Let $f:(\bc^2,0)\to(\bc,0)$ 
be an irreducible plane curve 
singularity with Newton pairs $\{(p_i,q_i)\}_{i=1}^s$ and let $n$ be an
integer $\geq 2$. Let $\G^{min}(\xfn)$ be the minimal (good) resolution
graph of the hypersurface singularity $(\xfn,0):=(\{f(x,y)+z^n=0\},0)$.
Then the following facts hold:

(a) The integer $s$ is uniquely determined by  $\G^{min}(\xfn)$.

(b) $s=1$ 
if and only if  $\G^{min}(\xfn)$ is either a string (with all the genera
 zero), or a star-shaped graph (where only the central vertex
might have genus $\g$ non-zero). Moreover,   $f(x,y)+z^n$
has the same equisingularity type as the Brieskorn singularity 
$x^{q_1}+y^{p_1}+z^n$.

$\G^{min}(\xfn)$ is a string if and only if $\{q_1,p_1,n\}=\{q_1,2,2\}$. 
If $\G^{min}(\xfn)$ is a star-shaped graph with $\g=0$, then the set of
integers $\{q_1,p_1,n\}$ is uniquely determined. Moreover, the 
only ambiguity which can appear in the case $\g>0$ is described in \ref{S1}.

(c) If $s=2$ then it can happen that two pairs $(f_1,n_1)$ and $(f_2,n_2)$
(but not more) 
provide identical graphs $\G^{min}(\xfn)$. If this is the case then both
 of them should satisfy the numerical restrictions: 
\begin{equation*}
(n,a_2)=1, \ \ \mbox{and}\ \ \  (n,p_2)a_1p_1|n.
\tag{$*$}
\end{equation*}
(which can be recognized from the graph as well), and 
(automatically) at least one of the vertices   has genus $\g>0$. 
This case is described in \ref{S2}.

(d) In all other cases (i.e.  for any $s\geq 3$ or for $s=2$ excluding the 
exceptional case ($*$)),
$\G^{min}(\xfn)$ determines uniquely  the Newton pairs
of $f$ and the integer $n$ (by a precise algorithm which basically constitutes
the next proof). 

(e) In particular, except for  the two cases S1  and S2 (cf. \ref{S1} and 
\ref{S2}), from the link one can recover completely the Newton pairs
of $f$ and the integer $n$ (provided that  we disregard the 
$z$-axis ambiguity, cf. \ref{zamb}.) In particular, this is true without any 
exception provided that the link is a rational homology sphere.

On the other hand, in the cases S1 and S2, the link together with 
the Milnor number of the hypersurface singularity $f+z^n$ determines 
completely the Newton pairs of $f$ and the integer $n$ (cf.
the two addendums in \ref{S1} and \ref{S2}). }

\begin{proof}
We denote  $\G^{min}(\xfn)$ by $\G$ and its rupture vertices 
by $\calr$. It is convenient to separate those cases when $\#\calr$ is small.

\vspace{2mm}

\noindent {\bf (Case A)} \ {\em Assume that $\calr=\varnothing$.} \ 
By \ref{ru} the set of rupture vertices  of $\Gm$ is never empty. 
Hence, by \ref{nonmin}, $\calr=\varnothing$ if and only if in the ``P-case''
we contract $v'$, and $v'$ is the unique rupture vertex of $\Gm$.  But
$\Gm$ has a unique rupture point if and only if $s=1$. 
Therefore (cf. \ref{nonmin}) this situation occurs if and only if 
$s=1, \ p_1=n=2$ and $(q_1,2)=1$; i.e. $f(x,y)+z^n$ has the equisingularity 
type of $x^{q_1}+y^2+z^2$. 
Clearly,  $q_1$  can be recovered from the graph: it is 
its determinant. Cf. also with \ref{S1}. 

\vspace{2mm}

\noindent {\bf (Case B)}\ {\em Assume that $\#\calr=1$.} \
From \ref{nonmin} is clear that in the ``P-case'' the number of rupture 
vertices of $\G$ is even. Hence, this case is
 excluded, and  the number of rupture 
vertices of $\Gm$ is also 1. By \ref{ru}, this can happen only of $s=1$.
In particular, $f+z^n$ is of Brieskorn type: $x^{q_1}+y^{p_1}+z^n$,
with $q_1>p_1\geq 2, \ (p_1,q_1)=1$ and $n\geq 2$ (where the case 
$p_1=n=2$ is excluded, see above). This case is completely covered by \ref{S1}.

\vspace{2mm}

\noindent {\bf (Case C)}\ {\em Assume that $\#\calr>1$.} By \ref{ru} and 
\ref{nonmin} $\#\calr>1$ if and only if $s\geq 2$. The proof (algorithm)
consists of several steps, each step recovers some data.

\vspace{2mm}

\noindent {\bf (1)} {\em The set $\calr_1$} can be determined from 
\ref{lgr} and \ref{iandii}. Indeed, we start with the set $\call(\G(\calr))$
(where $\G=\G^{min}(\xfn)$). Then, if there exists at least one $v\in 
\call(\G(\calr))$ with $\g_v\not=0$, then $\calr_1=\{v\in
\call(\G(\calr))\,:\, \g_v\not=0\}$ (cf. \ref{iandii}(a)).
If $\g_v=0$ for all $v$, then we consider their $\beta$-invariants $\beta(v)$.
If they are all equal, then by \ref{iandii}(b) one gets $\calr_1=
\call(\G(\calr))$. If they are not all equal, then only one can be $\leq 1/2$
(corresponding to $\{v\}=\calr_s$), and all the others are $>2$ (and equal to
 each other) corresponding to $\calr_1$ (cf. \ref{iandii}(b)). 

\vspace{2mm}

\noindent {\bf (2)}
 {\em The sets $\pi(q^{-1}(v_k))$ ($1\leq k\leq s$).} \ Define a 
{\em distance}
on the set $\calr$. If $w_1,\ w_2\in \calr$, and on the shortest 
path  in $\G$ connecting them there are exactly $l$ rupture vertices of $\G$
(including $w_1$ and $w_2$),
then we say that $d(w_1,w_2):=l-1\geq 0$. 
Moreover, for any subset $\calr'\subset \calr$ and $w\in\calr$ we define 
$d(\calr',w)$ as usual by $\min\{\, d(w',w)\,:\, w'\in\calr'\}$.

Then, for any $k\geq 1$, we write $\calr_k:=\{v\in \calr\, :\, 
d(\calr_1,v)=k-1\}$. Let $s':=\max\{k\,:\, \calr_k\not=\varnothing\}$. 
We distinguish two cases:

(a) \ $\#\calr_{s'}=1$. Then $s=s'$ and $\pi(q^{-1}(v_k))=\calr_k$ 
for $1\leq k\leq s$.

(b) \ $\#\calr_{s'}>1$. This can happen exactly in the ``P-case'' (cf.
\ref{nonmin}). In this situation, $s=s'+1$, $\pi(q^{-1}(v_k))=\calr_k$ for
$1\leq k\leq s-1$, and the (unique) vertex $v'= q^{-1}(v_s)$  of
$\Gm$  is ``missing''
in $\G$, i.e. $\pi(\{v'\})=\varnothing$ (cf. \ref{nonmin}). 

For a moment we postpone the ``P-case'', and we assume (a). We will come back
to the ``P-case'' in (11).

\vspace{2mm}

\noindent {\bf (3)} {\em The integers $\{h_k\}_{k=2}^s$} are determined by the 
identities $h_k=\#\calr_{k-1}/\#\calr_k$. ($h_1$ will be 
determined later.)

\vspace{2mm}

\noindent {\bf (4)} {\em The sets $\{\pi(q^{-1}(\bar{v}_k))\}_{k=2}^{s-1}$
and the integers $\{\tilde{h}_k\}_{k=2}^{s-1}$, \ $\{p_k\}_{k=2}^{s-1}$ 
(for $s\geq 3$)}.  Fix $2\leq k\leq s-1$ and some $v\in \calr_k$.
Consider the set of strings $\calst(v)$ supported by $v$ (cf. \ref{defs}(3)).

If $\calst(v)\not=\varnothing$, then $\#\calst(v)=\tilde{h}_k$ and 
$\det:\calst(v)\to\bfz$ is constant with value $p_k'$. Then $p_k=p_k'
\cdot h_k$. Then this is happening for any choice of $v$, and 
 $\pi(q^{-1}(\bar{v}_k))$ is the set of leaf vertices  of $\G$
situated on the strings of type $\cup_v\calst(v)$, $v\in \calr_k$. 

If $\calst(v)=\varnothing$ then all the strings of type $\G(\bar{v}_k)$
are collapsed in $\G$, in particular $\pi(q^{-1}(\bar{v}_k))=\varnothing$.
Hence their determinants $p_k'=p_k/h_k=1$.
In particular, $p_k=h_k\geq 2$. Then $\tilde{h}_k$ is given by the genus
formula $\tilde{h}_k-1=2\g_v/(h_k-1)$.

\vspace{2mm}

For the ``ends'' $k=1$ and $k=s$ we need more special computations
(since we have to separate the two different types of strings (which may be,
or may not  be ``missing'' from $\G$). In (5) we recover $\tilde{h}_s$, in
(6)  and (7) $n$ and $p_s$ and the arrowhead of $\Gm$ (excepting the case
($*$)). In (8) we treat the invariants with index $k=1$. 

\vspace{2mm}

\noindent {\bf (5)} {\em The integer $\tilde{h}_s$}. If $h_s>1$ then 
the genus formula for $\g_{v_s}$ gives $\tilde{h}_s$. If $h_s=1$, the strings
 of type  $\G(\bar{v}_s)$ cannot be collapsed, hence $\calst(v)\not
=\varnothing$ for $\{v\}=\calr_s$. Then $\tilde{h}_s=\alpha(v)$, 
cf. \ref{defs}(3) and \ref{c2}(b).

\vspace{2mm}

\noindent {\bf (6)}
 {\em $p_s$ and $n$ and the arrowhead of $\Gm$ in the following cases:}
$$\begin{array}{l}
\mbox{(i) \ \ either $s\geq 3$, or }\\
\mbox{(ii) \ \ $s=2$ but $\{a_1',p_1'\}\not=\{1\}$}\end{array}$$
First we show that in both cases we can compute the product $h_{s-1}\tilde{h}
_{s-1}$. Indeed, in the case (i) this follows from (3) and (4). If $s=2$ then
  we  proceed as follows. Since $a_1'$ and $p_1'$ are not both 1,
$\calst(v)\not=\varnothing$ for $v\in \calr_1$. If the determinant
has two values on this set, then $h_1\tilde{h}_1=\alpha(v)$, cf.
\ref{defs}(3) and \ref{c2}(a). If all the determinants are the same, then 
either $\G(\bar{v}_0)$ or  $\G(\bar{v}_1)$ is collapsed. If $\G(\bar{v}_0)$
is collapsed then $a_1'=1$ hence $\tilde{h}_1\geq 2$. In the second case
$h_1\geq 2$. Hence in both cases the following procedure works: 
take $c_1:=\#\calst(v)$ (which automatically is $\geq 2$), compute 
$c_2$ by the genus formula $2\g_{v_1}=(c_1-1)(c_2-1)$ and set
$h_1\tilde{h}_1=c_1c_2$. 

Now, we go back to $p_s$ and $n$ and the position of the arrowhead.

Notice that by \ref{ru} and part (2), the determinants of type
$D(v_{s-1})$ and $D_+(v_{s-1})$ are well-defined in $\G$, and their values do 
not change by the minimalization procedure. E.g. 
$D_+(v_{s-1})$ can be computed from \ref{dd}(b).
Notice also that $D_{St}(v_s)=n(p_s')^{\tilde{h}_s}/h_s\tilde{h}_s$ (cf. 
\ref{pr}).
Therefore
\begin{equation*}
\frac{D(v_{s-1})^{h_s-1}D_{St}(v_s)}{D_+(v_{s-1})}=\frac{h_{s-1}\tilde{h}_{
s-1}p_s}{h_s\tilde{h}_s}.
\end{equation*}
Hence this value can be determined from the graph, a  fact which is true for
$h_{s-1}\tilde{h}_{s-1}$ (see above) and $h_s$ (cf. (3)) and $\tilde{h}_s$ 
(cf. (5)) as 
well. Hence we get $p_s$. In particular, we can compute the string 
determinants $D(\bar{v}_s)=p_s'$ and (using $D_{St}(v_s)$) $D(v_s)=n/(h_s\tilde{h}_s)$
as well. This gives $n$ too. If $D(v_s)\not=1$ then we put the arrow on the 
string with this determinant (cf. \ref{c2}); if $D(v_s)=1$ then we put the 
arrowhead on  $\pi(q^{-1}(v_s))$. In this way we recover the arrow of $\Gm$.

\vspace{2mm}

\noindent {\bf (7)} {\em $p_s$ and $n$ and the arrowhead of $\Gm$ if }
$$\mbox{(iii)} \ \ \ \tilde{h}_s\not=1.$$
Consider  $\calst(v)$ for $\{v\}=\calr_s$, cf. \ref{defs}(3). 
Then by a verification $\cald_{St}^{red}(v)=
(p_s')^{\tilde{h}_s-1}$. Since $\tilde{h}_s$ is determined in (5), and 
it is $\not=1$, one gets $p_s'$. Then we repeat the arguments of 
(6).

\vspace{2mm}

\noindent {\bf (8)} {\em The integers $a_1, \ p_1, \ h_1$ and $\tilde{h}_1$ in 
the cases when one of the conditions (i) or (ii)  or (iii) is valid.} 
Fix a vertex $v\in\calr_1$ and consider $\calst(v)$ as in \ref{defs}(3).

If $\calst(v)\not=\varnothing$, then $\#I(v)\leq 2$. 
If $\#I(v)=2$, then compute
the two numbers $D_1\cdot \#_2$ and $D_2\cdot \#_1$.
They are the candidates for $a_1=q_1$ and $p_1$, cf. \ref{pr}. 
Since $q_1>p_1$, these two numbers cannot be the same. If, say, 
$D_1\cdot \#_2>D_2\cdot \#_1$, then 
$\calst_1$ is the index set of $\pi(q^{-1}(\bar{v}_0))$ and 
$\calst_2$ of $\pi(q^{-1}(\bar{v}_1))$. Hence $h_1=\#_1$,
$\tilde{h}_1=\#_2$, $q_1=D_1\cdot \#_2$ and $p_1=D_2\cdot \#_1$. 

If there is only one level set with data $D_1$ and $\#_1$, then in the above 
argument we write $D_2=1$ and we determine $\#_2$ using the genus formula
$2\g_{v_1}=(\#_1-1)(\#_2-1)$ (which is possible since $D_2\#_1=\#_1\geq 2$). 
And we repeat the above argument. 

Now we assume that $\calst(v)=\varnothing$. This can happen only if 
$a_1'=p_1'=1$, hence $q_1=\tilde{h}_1$ and $p_1=h_1$. 
First we determine $H:=h_1\tilde{h}_1$.

$D(v_1)$ gives an equation of type $q_2=H\cdot A$, where $A$ is a positive
number which can be determined from the graph by the previous steps.
Moreover, $D_-(v_2)=a_2'$, hence $a_2=D_-(v_2)\tilde{h}_2$ is known from the 
graph. Finally, $a_1p_1p_2=Hp_2$, where $p_2$ too is known from the graph.
Then the identity $a_2=q_2+a_1p_1p_2$ gives a non-trivial linear
equation for $H$. 

Then $h_1\tilde{h}_1=H$ and $(h_1-1)(\tilde{h}_1-1)=2\g_{v_1}$ provides 
$h_1$ and $\tilde{h}_1$ modulo their  permutation.  But 
$\tilde{h}_1=q_1>a_1=h_1$, hence we get $h_1$ and $\tilde{h}_1$. 

\vspace{2mm}

\noindent {\bf (9)} {\em The integers $\{a_k\}_{k=2}^s$  when
one of the conditions (i) or (ii)  or (iii) is valid.} 
Once we have the position of the arrow, 
we have all the multiplicities $\{\m_{v_k}\}_k$ (cf. \ref{1.3}), hence 
\ref{l2}(b) gives all the integers $a_k'$. An alternative way is to use 
inductively \ref{gen}. 

\vspace{2mm}

\noindent {\bf (10)}
 {\em Assume that the conditions (i), (ii) and (iii) are not
valid.} This means that $s=2$ and $a_1'=p_1'=\tilde{h}_2=1$.
This is exactly the case of S2-coincidence treated in \ref{S2}.

\vspace{2mm}

\noindent {\bf (11)}
 {\em The ``P-case''.} Now we go back to step (2), case (b).
In this case $\#\calr_{s'}=2$, so write  $\calr_{s'}=\{w_1,w_2\}$. 
Take the shortest path in $\G$ connecting $w_1$ and $w_2$. Take the edge 
 ``at the  middle of the path'', blow it up, and put an arrow on it.
This new graph is exactly $\Gm$. Set $\calr_s:=\{v'\}$, where $v'$ is the
new  vertex. Then we can repeat all the above arguments. 

Notice that step (6) works, since if $s=2$ and $a_1'=p_1'=1$, then 
$h_1=p_1\geq 2$ and $\tilde{h}_1=a_1\geq 2$ hence $\g_{v_1}>0$.
But in the ``P-case'' all the genera are zero. (Hence step (7) is not needed.)
[In fact, since in this case we already have the position of the arrow, 
we can compute some of the invariants much faster using the multiplicities and
\ref{l2}(b).] \end{proof}

\subsection{}\label{sim1}  In the above proof we were rather meticulous in
separating the possible sets $\pi(q^{-1}(v))$. The fruit of this is the 
following corollary (whose  proof is left to the reader, and  basically it is
incorporated in the previous proof of the main theorem). 

First recall that the cyclic covering $\xfn\to X$ has a $\bz_n$ Galois action.
This lifts to the level of the resolution, hence $\G^{min}(\xfn)$ inherits
a natural $\bz_n$-action as well. The question is: has the graph $\G^{min}
(\xfn)$ any extra symmetry? 

Take for example the Brieskorn case \ref{Seifert}. Then the Galois action 
permutes cyclically (via its image $\bz_h$)
the $h$ arms with Seifert invariants $a$. On the other 
hand, the total symmetry group of the graph is the total permutation 
group of these arms. So, in this sense, the symmetry group of the graph is
definitely larger than the (image) of the Galois action. On the other hand, 
their orbits are the same. This fact is valid in general.

\subsection{Corollary.}\label{sim2} {\em Assume that $\sigma$ is a 
(decorated graph-) automorphism of $\G^{min}(\xfn)$ which identifies 
two vertices, say, $v_1$ and $v_2$. Then $v_1$ and $v_2$ are in the same 
orbit with respect to the Galois action. }

\vspace{2mm}

This result definitely cannot be extended to the general case when $f$ is not
irreducible. E.g., $\G^{min}(\bc^2,0)$ can have a symmetry (take e.g.
$f=(x^2+y^3)(x^3+y^2)$) which lifts to an automorphism of $\G^{min}(\xfn)$,
which does not come from the Galois covering.

\end{document}